\documentclass[12pt]{article}
\usepackage{enumerate,eepic,epic,amssymb,amsfonts,times,color,epsfig,
fancybox,graphicx,amsmath,pifont,pxfonts,txfonts,bbm,mathdots,mathtools}
\makeatletter
\ifcase\@ptsize
 \font\teneufm=eufm10
 \font\seveneufm=eufm7
 \font\fiveeufm=eufm5
 \font\teneusm=eusm10
 \font\seveneusm=eusm7
 \font\fiveeusm=eusm5
\or
 \font\teneufm=eufm10 scaled \magstephalf
 \font\seveneufm=eufm7
 \font\fiveeufm=eufm5
 \font\teneusm=eusm10 scaled \magstephalf
 \font\seveneusm=eusm7
 \font\fiveeusm=eusm5
\or
 \font\teneufm=eufm10 scaled \magstep1
 \font\seveneufm=eufm7
 \font\fiveeufm=eufm5
 \font\teneusm=eusm10 scaled \magstep1
 \font\seveneusm=eusm7
 \font\fiveeusm=eusm5
\fi

\newfam\eufmfam
\newfam\eusmfam
\textfont\eufmfam=\teneufm  \scriptfont\eufmfam=\seveneufm
  \scriptscriptfont\eufmfam=\fiveeufm
\textfont\eusmfam=\teneusm  \scriptfont\eusmfam=\seveneusm
  \scriptscriptfont\eusmfam=\fiveeusm

\def\frak{\ifmmode\let\next\frak@\else
 \def\next{\errmessage{Use \string\frak\space only in math mode}}\fi\next}
\def\frak@#1{{\frak@@{#1}}}
\def\frak@@#1{\fam\eufmfam#1}

\def\sh{\ifmmode\let\next\sh@\else
 \def\next{\errmessage{Use \string\sh\space only in math mode}}\fi\next}
\def\sh@#1{{\sh@@{#1}}}
\def\sh@@#1{\fam\eusmfam#1}

\ifcase\@ptsize
 \font\tenmsa=msam10
 \font\sevenmsa=msam7
 \font\fivemsa=msam5
 \font\tenmsb=msbm10
 \font\sevenmsb=msbm7
 \font\fivemsb=msbm5
\or
 \font\tenmsa=msam10 scaled \magstephalf
 \font\sevenmsa=msam7
 \font\fivemsa=msam5
 \font\tenmsb=msbm10 scaled \magstephalf
 \font\sevenmsb=msbm7
 \font\fivemsb=msbm5
\or
 \font\tenmsa=msam10 scaled \magstep1
 \font\sevenmsa=msam7
 \font\fivemsa=msam5
 \font\tenmsb=msbm10 scaled \magstep1
 \font\sevenmsb=msbm7
 \font\fivemsb=msbm5
\fi

\newfam\msafam
\newfam\msbfam
\textfont\msafam=\tenmsa  \scriptfont\msafam=\sevenmsa
  \scriptscriptfont\msafam=\fivemsa
\textfont\msbfam=\tenmsb  \scriptfont\msbfam=\sevenmsb
  \scriptscriptfont\msbfam=\fivemsb

\def\Bbb{\ifmmode\let\next\Bbb@\else
 \def\next{\errmessage{Use \string\Bbb\space only in math mode}}\fi\next}
\def\Bbb@#1{{\Bbb@@{#1}}}
\def\Bbb@@#1{\fam\msbfam#1}
\def\hexnumber@#1{\ifnum#1<10 \number#1\else
 \ifnum#1=10 A\else\ifnum#1=11 B\else\ifnum#1=12 C\else
 \ifnum#1=13 D\else\ifnum#1=14 E\else\ifnum#1=15 F\fi\fi\fi\fi\fi\fi\fi}
\def\msa@{\hexnumber@\msafam}
\def\msb@{\hexnumber@\msbfam}
\mathchardef\square="0\msa@03

\makeatother
\newcommand{\beq}{\begin{equation}}
\newcommand{\eeq}{\end{equation}}
\newcommand{\ba}{\begin{array}}
\newcommand{\ea}{\end{array}}
\newcommand{\bea}{\begin{eqnarray}}
\newcommand{\eea}{\end{eqnarray}}
\newcommand{\bean}{\begin{eqnarray*}}
\newcommand{\eean}{\end{eqnarray*}}

\newtheorem{theorem}{Theorem}[section]
\newtheorem{prop}[theorem]{Proposition}
\newtheorem{lem}[theorem]{Lemma}
\newtheorem{defi}[theorem]{Definition}
\newtheorem{cor}[theorem]{Corollary}
\newtheorem{remark}[theorem]{Remark}

\newtheorem{proof}{Proof.}
\newenvironment{pf}{\begin{proof} \rm}{\hfill$\square$ \end{proof}}

\makeatletter
\@addtoreset{equation}{section}

\makeatother

\def\endpf{\begin{flushright}$\square$\end{flushright}}

\begin{document}
\begin{titlepage}
\begin{center}
{\huge The classification of the cyclic 
$\mathfrak{sl}(n+1)\ltimes \mathbbm{C}^{n+1}$--modules}
\end{center}
\vspace{0.8truecm}
\begin{center}
{\large
Paolo Casati
\vskip  0.8truecm
Dipartimento di Matematica e applicazioni\\
 Universit\`a di Milano-Bicocca\\
Via Cozzi 53, I-20125 Milano, Italy\\}
\end{center}
E--mail:  paolo.casati@unimib.it\\
\vspace{0.2truecm}
\vspace{0.2truecm}
\abstract{In this paper we classify all the cyclic finite dimensional 
indecomposable\\ modules of the perfect Lie algebras $\mathfrak{sl}(n+1)\ltimes 
\mathbbm{C}^{n+1}$, given by the semidirect sum of the simple Lie algebra $A_n$  
with its standard representation. 
Furthermore, using the embeddings of the Lie algebras  $\mathfrak{sl}(n+1)\ltimes 
\mathbbm{C}^{n+1}$ in  $\mathfrak{sl}(n+2)$, 
 we show that any finite dimensional irreducible module of 
$\mathfrak{sl}(n+2)$ restricted to 
$\mathfrak{sl}(n+1)\ltimes \mathbbm{C}^{n+1}$ is a  cyclic module
and that any cyclic $\mathfrak{sl}(n+1)\ltimes \mathbbm{C}^{n+1}$--modules 
can be constructed as quotient module  of the restriction to  
$\mathfrak{sl}(n+1)\ltimes \mathbbm{C}^{n+1}$ of some finite dimensional 
irreducible $\mathfrak{sl}(n+2)$--modules. This explicit realization of the 
cyclic $\mathfrak{sl}(n+1)\ltimes \mathbbm{C}^{n+1}$--modules plays a role in their classification. } 

\end{titlepage} 
\section{Introduction}
Surely the  construction  and the (even partial) classification of  the 
indecomposable modules of a given Lie algebra is a very hard and (in most of the cases) an almost unsolved problem. In fact only for   the class of all semisimple Lie algebras such classification has been accomplished, while, for the non  semisimple Lie algebras very few results have  been, until now,  achieved.\par 
In the case of semisimple Lie algebra the classification of all indecomposable 
module  can be summarized by saying that any  indecomposable module is 
an irreducible modules and for a given simple Lie algebra $\mathfrak{s}$ 
the set of irreducible modules is classified by an element in $\mathbbm{N}^n$,
where  $\mathbbm{N}$ is the set of the  natural number and $n$ 
is the rank of $\mathfrak{s}$. In fact this element of  $\mathbbm{N}^n$
is a dominant weight of $\mathfrak{g}$ see \cite{H}, \cite{FH}  or the next section for more 
details.\par 
In order to obtain similar result for non semisimple Lie algebra one has  
 to identify a distinguished class of indecomposable
representations for which one could expect to obtain a reasonable 
classification.  This seems to be possible in two ways.\par 
First, one can   consider the embeddings of a given Lie algebra 
$\mathfrak{g}$  into a semisimple Lie algebra $\mathfrak{s}$  
(or other Lie algebras whose indecomposable representations are 
at least partial known, see for instance \cite{CMS} for an example where a 
not semisimple Lie algebra, in fact a truncated current Lie algebra \cite{CO}  is used)
and to use the well known irreducible finite dimensional modules of $\mathfrak{s}$  to study and    classify the indecomposable $\mathfrak{g}$--modules obtained by restriction. 
Such approach   to the theory of the indecomposable module of non semisimple Lie algebra has been pioneered  by Douglas and Premat in \cite{DP} and in the PhD Thesis of Douglas and further developed by  Deguise, Douglas,  Premat, Repka  
 in  \cite{D},\cite{DdG},\cite{P},\cite{DR} and by Minniti, Salari and the author in \cite{CMS}. \par 
Second, one can select a particular class of indecomposable modules
which satisfy suitable properties like to be cyclic or uniserial.
 In this contest in fact Cagliero and Szechtman \cite{CF} in a recent beautiful 
paper have classified   all uniserial $\mathfrak{g}$--modules
of  the perfect Lie algebras  $\mathfrak{g} = \mathfrak{sl}(2)\ltimes V(m)$, 
where $V(m)$ is the irreducible 
$\mathfrak{sl}(2)$--module with highest weight $m\geq  1$,
Previously only the case corresponding to the Lie algebra  
$\mathfrak{sl}(2)\ltimes V(1)$,  \cite{Pi}
(where actually all the indecomposable representations with one generator had 
been classified), was known.\par 
The aim of this paper is to proceed further into these lines of though, using 
both approaches. 
We shall indeed consider  the perfect Lie algebras 
$\mathfrak{sl}(n+1)\ltimes \mathbbm{C}^{n+1}$  given by the semidirect sum of 
the simple Lie algebra $A_n$  with its standard representation.\par 
The main results  is a complete classification of
all cyclic finite dimensional 
$\mathfrak{sl}(n+1)\ltimes \mathbbm{C}^{n+1}$--module with $n\geq 1$
(observe that the case $n=1$ was already obtained by Piard in \cite{Pi}),
i.e., the indecomposable modules with only one generators, 
and therefore as subcase  all the uniserial 
$\mathfrak{sl}(n+1)\ltimes \mathbbm{C}^{n+1}$--modules, see Theorem  
\ref{classificationth}. We shall namely show that such indecomposable modules
are classified  by particular bounded subsets of integer numbers see   
Definition \ref{setM}. \par 
In order to achieve such result we shall make use of the fact that 
the perfect Lie algebra  $\mathfrak{sl}(n+1)\ltimes \mathbbm{C}^{n+1}$
can be embedded in the simple Lie algebra $\mathfrak{sl}(n+2)$ see
\cite{DP}. \par 
More precisely we first show how to associate a set $\mathcal{J}(V)$  of integer numbers 
to any cyclic  $\mathfrak{sl}(n+1)\ltimes \mathbbm{C}^{n+1}$--module $V$ and that this correspondence is injective: to   inequivalent cyclic modules correspond  different subsets of integer numbers see Theorems \ref{m1mn} and \ref{classificationp1}. 
Then we use the restriction to 
$\mathfrak{sl}(n+1)\ltimes \mathbbm{C}^{n+1}$ of the finite dimensional 
irreducible $\mathfrak{sl}(n+2)$--modules in order to construct 
for any such set of integer numbers  a corresponding 
cyclic $\mathfrak{sl}(n+1)\ltimes \mathbbm{C}^{n+1}$--module see Theorem 
\ref{setjlambda}.
Here it is interesting to note that any 
$\mathfrak{sl}(n+1)\ltimes \mathbbm{C}^{n+1}$ module obtained by 
restriction of a finite dimensional 
irreducible $\mathfrak{sl}(n+2)$--modules $V$  on 
$\mathfrak{sl}(n+1)\ltimes \mathbbm{C}^{n+1}$ is itself a 
cyclic module generated by the highest weight vector of $V$.\par 
The paper is organized as follows.
In Section 2 we collect all the needed properties of the perfect Lie algebras 
$\mathfrak{sl}(n+1)\ltimes \mathbbm{C}^{n+1}$ and of their embeddings in 
$\mathfrak{sl}(n+2)$. Furthermore we  brifely describe the  finite irreducible 
module of the simple Lie algebra of type $A_n$ together with the basis of such 
modules recently founded by Feigin,  Fourier and Littelmann \cite{FFL}, which
will play a important role in the last section.
In Sections 3  we first present the definition of cyclic module together with some general result on cyclic 
module of perfect Lie algebras. Then we obtain a detailed description of the 
finite dimensional cyclic modules of the perfect Lie algebra  
$\mathfrak{sl}(n+1)\ltimes \mathbbm{C}^{n+1}$ as direct sum of irreducible 
$\mathfrak{sl}(n+1)$--modules  and determine how the  radical of 
$\mathfrak{sl}(n+1)\ltimes \mathbbm{C}^{n+1}$ intertwine them. 
Finally we show how any  finite dimensional cyclic module $V$  is determined by a bounded subset of integer numbers $\mathcal{J}(V)$. In Section 4 we study the 
indecomposable modules of $\mathfrak{sl}(n+1)\ltimes \mathbbm{C}^{n+1}$
which arise a restriction to them of the finite irreducible 
$\mathfrak{sl}(n+2)$--modules,  showing that they are cyclic and that any 
cyclic $\mathfrak{sl}(n+1)\ltimes \mathbbm{C}^{n+1}$--module can be viewed as 
their quotient. This latter result together with those obtained in section 3 
complete the classification of the $\mathfrak{sl}(n+1)\ltimes 
\mathbbm{C}^{n+1}$--modules.
 \section{The perfect Lie algebra $\mathfrak{sl}(n+1)\ltimes \mathbbm{C}^{n+1}$} 
In this section we will describe into the details needed for our purposes the 
perfect Lie algebra 
$\mathfrak{sl}(n+1)\ltimes \mathbbm{C}^{n+1}$ given  by the semidirect sum of 
the classical simple Lie algebra $\mathfrak{sl}(n+1)$  with its standard 
representation \cite{DR1}.\par \noindent  
Throughout this paper, all the Lie algebras and their modules (which are supposed to be always finite dimensional)  are over the 
ground field  $\mathbbm{C}$ of the complex numbers.\par
Let $\mathfrak{s}$ be  a simple complex Lie algebra, $\mathfrak{h}\subset\mathfrak{s}$  a Cartan subalgebra,  $\mathfrak{h}^*$ its 
complex dual
and $\Delta=\Delta(\mathfrak{s},\mathfrak{h})\subset \mathfrak{h}^*$ the  set of 
roots of $\mathfrak{h}$ in $\mathfrak{s}$. 
Denote by $(\ \  ,\ )$ the Killing form on $\mathfrak{g}$. The induced form on $\mathfrak{h}^*$  will be denoted by $(\ \  ,\ )$ as well.\par 
For   $\alpha \in \Delta$, let $H_\alpha\in \mathfrak{h}$ be the corresponding coroot and let 
  $\mathfrak{s}_\alpha=\{X\in \mathfrak{s}\ \vert \ \left[ H,X\right]=\alpha(H)X\ \forall H\in  \mathfrak{h}\}$   its  root--space.\par
A set of simple roots $\Pi=\left\{\alpha_1,\dots,\alpha_n\right\}\subset\Delta$  is a basis of $\mathfrak{h}^*$  such that any root $\alpha$ belongs either to the  
monoid $\mathbbm{Z}_{\geq 0}\Pi$ or to its opposite. Accordingly we can view $\Delta$ as  the disjoint union of positive roots $\Delta_\Pi^+=\mathbbm{Z}_{\geq 0}\Pi\cap \Delta$ 
and negative roots $\Delta_\Pi^-=-\Delta_\Pi^+$. (Since we are going to choose a set of simple roots once for ever we shall simply write $\Delta^\pm$ instead of 
$\Delta_\Pi^\pm$).
 This  decomposition of $\Delta$ induces the Cartan decomposition of $\mathfrak{s}$:
$$
\mathfrak{s}=\mathfrak{n}^+\oplus\mathfrak{h}\oplus\mathfrak{n}^-
$$
where $\mathfrak{n}^+=\oplus_{\alpha\in \Delta^+}\mathfrak{s}_\alpha$ and   $\mathfrak{n}^-=\oplus_{\alpha\in \Delta^-}\mathfrak{s}_\alpha$ are maximal nilpotent subalgebras.\par\noindent
 Let $E_\alpha$  be a basis of $\mathfrak{s}_\alpha$ and let $F_\alpha$  in  
$\mathfrak{s}_{-\alpha}$  be defined by the requirement $\left[E_\alpha,F_\alpha\right]= H_\alpha$. \par\noindent 
Denote by $\mathcal{U}(\mathfrak{s})$, $\mathcal{U}(\mathfrak{n}^+)$, $\mathcal{U}(\mathfrak{n}^-)$ the universal enveloping algebras of 
$\mathfrak{s}$, $\mathfrak{n}^+$, $\mathfrak{n}^-$ respectively. (More in general $\mathcal{U}(\mathfrak{g})$ will denote the universal enveloping algebra of the  Lie algebra $\mathfrak{g}$.)\par
If the simple Lie algebra $\mathfrak{s}$ is the complex Lie algebra $\mathfrak{sl}(n+1)$, we can choose as Cartan subalgebra the set of all diagonal matrices in 
$\mathfrak{sl}(n+1)$, then the Cartan decomposition becomes the usual  triangular decomposition of  $\mathfrak{sl}(n+1)$ 
in into the direct sum of strictly upper triangular, diagonal, and strictly lower triangular matrices.\par
If $H$ is the diagonal matrix with $h_1,\dots ,h_{n+1}$ on the main diagonal we set $\epsilon_j(H) =h_j$. Then 
$\Delta(\mathfrak{g},\mathfrak{h})=\{\epsilon_i-\epsilon_j,\ \ i,j=1,\dots n+1,\ \ \   i \neq j\  \}$ and we may take as set of simple roots 
$\Pi=\left\{\ \alpha_1=\epsilon_1-\epsilon_2,\ \alpha_2=\epsilon_2-\epsilon_3,\dots,\ \alpha_n= \epsilon_{n} -\epsilon_{n+1}\right\}$,  hence  the  roots of 
$\mathfrak{sl}(n+1,\mathbbm{C})$  are of the form $\pm(\alpha_p + \alpha_{p+1} +\dots  + \alpha_q)$ for some $1\leq p \leq q \leq n$. Following 
\cite{FFL} we  set for $1\leq p < q\leq n$:
$$
\alpha_{p,q}=\alpha_p+\dots \alpha_q\qquad H_{p,q}=H_{\alpha_{p,q}}\qquad E_{p,q}=E_{\alpha_{p,q}}\qquad F_{p,q}=F_{\alpha_{p,q}},
$$
and for our convenience also  $\alpha_{i,i} =\alpha_{i}$, $H_{i,i}=H_{\alpha_i}=H_i$, $E_{i,i}=E_{\alpha_i}=E_i$, and   $F_{i,i}=F_{\alpha_i}=F_i$.\par
Let   $A_{i,j}$ be   the $(n+1)\times (n+1)$ matrix with $1$ in the $i,j$ position and $0$'s everywhere else.  
Then  the coroot $H_{{p,q}}$ is the matrix $H_{{p,q}}=A_{p,p}-A_{q,q}$,  and the elements 
$E_{{p,q}}$ and $F_{{p,q}}$ are respectively  $E_{{p,q}}=A_{p,q+1}$ and $F_{{p,q}}=A_{q+1,p}$, $1\leq p\leq q\leq n$.\par
An element of $\mathfrak{h}^*$ is called a weight.
The set $P=\{\lambda\in  \mathfrak{h}^*\vert \ \lambda(H_\alpha)\in \Bbb Z,\ \forall \alpha \in \Delta\}$ is said the set of  integral
weights of $\mathfrak{s}$. A weight $\lambda$  of $P$ is said dominant if 
$\lambda(H_\alpha)\geq 0$ for any simple root $\alpha$, let us denote by 
$\Lambda$ the subset of the integral dominant weights.\par
A basis of $\mathfrak{h}^*$  is given by the fundamental weights 
$\left\{\omega_1,\dots,\omega_n\right\}$  defined by
the relations $\omega_i(H_j) = \delta_{i.j}$ for $i,j = 1\dots  n$,
where $\delta_{i.j}$ is the usual Kronecker delta.
With respect to this basis 
the integral dominant weights in $\mathfrak{h}^*$ can be written as 
$\lambda=\sum_{i=1}^n\lambda_i\omega_i$, $\lambda_i\in \mathbbm{N}$. 
It is also convenient to  set $\omega_{n+1}=\omega_0=0$.\par
The complex finite dimensional irreducible representations of 
$\mathfrak{sl}(n+1)$ are\\  parametrized by the dominant integral 
weights. We denote by   $V(\lambda)$ the   finite dimensional irreducible 
$\mathfrak{sl}(n+1)$--module corresponding to the    integral 
dominant weight $\lambda$, and by  $v_\lambda$ its    highest weight vector  i.e., a non trivial vector in $V(\lambda)$ of 
weight $\lambda$ annihilated by $\mathfrak{n}^+$: $Hv_\lambda=\lambda(H)v_\lambda$ for all $H\in \mathfrak{h}$ and  $\mathfrak{n}^+v_\lambda=0$,  which  generates $V(\lambda)$: 
$V(\lambda)=\mathcal{U}(\mathfrak{s})v_\lambda$.\par
An element $\mu$ of $\mathfrak{h}^*$ is said a weight of an  irreducible finite 
dimensional module $V(\lambda)$ if the weight space $V_\mu=\{v\in V(\lambda) 
\vert \ Hv=\mu(H)v\ \forall H\in \mathfrak{h}\}$ is different from the zero vector. Denote 
by $P(\lambda)$ the set of all weights of $V(\lambda)$. The module  $V(\lambda)$ 
may be decomposed as the direct sum of its weight spaces:
\beq\label{weigdec}
V(\lambda)=\bigoplus_{\mu\in P(\lambda)}V_\mu.
\eeq
In what follows we shall need the  explicit basis of $V(\lambda)$ (for the 
simple Lie algebra of type $A$.) constructed in a beautiful paper of Feigin, 
Fourier and  Littelmann \cite{FFL}.
This basis  conjectured by Vinberg is related to the notion of Dyck path.
\begin{defi}\label{Dyck} \cite{FFL}. A Dyck path  is a sequence
$$
\mathbf{p} = \left(\beta(0),\beta(1), \dots ,\beta(k)\right), \qquad  k\geq  0;
$$
of positive roots satisfying the following conditions:
\begin{enumerate}
\item  If $k = 0$, then $\mathbf{p}$ is of the form $\mathbf{p} = (\alpha_i)$ 
for some simple root $\alpha_i$.
\item If $k\geq 1$, then:
\begin{enumerate}
\item  the first  and last elements are simple roots. More precisely, $\beta(0) 
=\alpha_i$
and $\beta_k =\alpha_j$ for some $1\leq i < j\leq  n$.
\item  the elements in between obey the following recursion rule. If $\beta(s) 
=\alpha_{p,q}$
then the next element in the sequence is of the form either $\beta(s + 1) 
=\alpha_{p,q+1}$ or 
$\beta(s + 1) =\alpha_{p+1,q}$.
\end{enumerate}
\end{enumerate}
\end{defi}
 Let $\mathcal{S}(\mathfrak{n}^-)$  denote the symmetric algebra of 
$\mathfrak{n}^-$. Then for  
a multi--exponent $\mathbf{s} = (s_\beta)_{\beta\in \Delta^+}$, $s_\beta \in 
\mathbbm{Z}_{\geq 0}$, let $F^{\mathbf{s}}$
 be the element
$$
F^{\mathbf{s}}=\prod_{\beta \in 
\Delta^+}F_\beta^{s_\beta}\in\mathcal{S}(\mathfrak{n}^-).
$$
\begin{theorem}\label{dyckbasis}\cite{FFL}  Let $\lambda=\sum_{i=1}^n 
\lambda_i\omega_i$ be  an integral dominant $\mathfrak{sl}(n+1)$--weight 
and let $S(\lambda)$ be the set of all multi-exponents 
$\mathbf{s}=(s_\beta)_{\beta\in \Delta^+}\in \mathbbm{Z}^{\Delta^+}_{\geq 0}$  
such that, for all Dyck paths
$\mathbf{p}=(\beta(0),\dots \beta(k))$:
\beq
\label{dyckformula}
s_{\beta(0)}+s_{\beta(1)}+\dots+s_{\beta(k)}\leq \lambda_i+\lambda_{i+1}+\dots 
\lambda_j
\eeq
where $\beta(0) =\alpha_i$ and $\beta(k) =\alpha_j$. Then if   
 $v_\lambda$ is a highest vector of  $V(\lambda)$, the set 
$F^{\mathbf{s}}v_\lambda$  with $\mathbf{s}\in S(\lambda)$ forms a basis of $V(\lambda)$, which we shall call  the Feigin Frenkel Littelmann basis ((FFL) basis) of $V(\lambda)$.
\end{theorem} 
The standard module of $\mathfrak{sl}(n+1)$ on $\mathbbm{C}^{n+1}$ coincides with the highest weight  module $V(\omega_1)$. A suitable basis for this space 
together  with the action on it of the elements of  $\mathfrak{sl}(n+1)$
is  given by
\begin{prop}\label{omega1basis} 
Let $u_{\omega_1}$ be a highest weight vector of  $V(\omega_1)$ then 
\par\noindent 
1) the set:
\begin{equation}
\label{basis}
\mathcal{S}=\left\{v_{\omega_1}, F_1v_{\omega_1},\cdots,   F_{1i}v_{\omega_1},\cdots, F_{1n}v_{\omega_1}\right\}
\end{equation}
is a   basis  of $V(\omega_1)$.\par \noindent
2) The action of $\mathfrak{sl}(n+1)$ on $\mathcal{S}$  is given by the 
relations
\begin{equation}
\label{oegaaction}
\begin{array}{lll}
&H_jv_{\omega_1}=\delta_{j1}v_{\omega_1}\qquad  &H_jF_{1,i}v_{\omega_1}=
\delta_{j,i+1}F_{1,i}v_{\omega_1}-\delta_{j,i}F_{1,i}v_{\omega_1}\\
&E_{j-h,j}v_{\omega_1}=0\qquad &E_{j-h,j}F_{1,i}v_{\omega_1}=\delta_{j,i}F_{1,j-h-1}v_{\omega_1}, \quad F_{1,0}=1\\
&F_{j-h,j}v_{\omega_1}=\delta_{1,j-h}F_{1,j}v_{\omega_1}\quad &F_{j-h,j}F_{1,i}v_{\omega_1}=\delta_{j-1,i+h}F_{1,j}v_{\omega_1}\quad i,j=1,\dots, n \ \    0\leq h <j.\\
\end{array}
\end{equation}
\end{prop}
{\bf Proof} 1)  That the  set $\mathcal{S}$  (\ref{basis}) is a  basis of 
the $\mathfrak{sl}(n+1)$--module   $V(\omega_1)$  can be checked using Theorem  \ref{dyckbasis} \cite{FFL}. However,  in this very simple case,  it can be also shown
 directly. Indeed since for any $1\leq i\leq n$ the subalgebra spanned by the 
elements $\left\{H_{1i},E_{1i},F_{1i}\right\}$ is isomorphic to the simple Lie 
algebra  $\mathfrak{sl}(2)$ and since 
$H_{1i}v_{\omega_1}=\omega_1(\alpha_{1i})v_{\omega_1}=v_{\omega_1}$, 
 from the theory of the 
representations of $\mathfrak{sl}(2)$ follows that  $F_{1,i}v_{\omega_1}\neq 0$.
Further, since the element $F_{1,i}v_{\omega_i}$ has weight 
$\omega_i-\alpha_{1,i}$ ,
two different elements of the set $\mathcal{S}$   belong to two different 
weight--spaces. 
To show that the set $\mathcal{S}$ is a basis of $V(\omega_1)$ it suffices 
therefore to show that its  number of elements   is  
equals the dimension $n+1$ of  $V(\omega_1)$, which is obvious. \par\noindent
2) The equations (\ref{oegaaction}) follow immediately from the commutation 
relations of $\mathfrak{sl}(n+1)$, which for the convenience of the 
reader are explicitly given in the next proposition. (see equation 
(\ref{comrel})).
\endpf 
Therefore setting $P_1=v_{\omega_1}$, $P_{j+1}=F_{1,j}v_{\omega_i}$, 
$1\leq j \leq n$ and using the basis of $\mathfrak{sl}(n+1)$ given above we have  
\begin{prop}\label{slncn}  The perfect Lie algebra $\mathfrak{sl}(n+1)\ltimes 
\mathbbm{C}^{n+1}$ is the Lie algebra spanned by the elements 
$$
\left\{ H_i,E_{p,q},F_{p,q}, P_{j}\right\}\qquad 1\leq i\leq n, \quad 1\leq p\leq q\leq n,\quad 1\leq j\leq n+1
$$
whose non trivial Lie brackets are:
\begin{equation}
\label{comrel}
\begin{array}{lll}
&\left[H_i,E_{p,q}\right]=(\alpha_i,\alpha_{p,q}) E_{p,q}&\quad
\left[H_i,F_{p,q}\right]=-(\alpha_i,\alpha_{p,q}) F_{p,q}\\
&\left[E_{p,q},E_{r,s}\right]=\delta_{q,r-1}E_{p,s}-\delta_{p,s+1}E_{r,q} &\quad 
 \left[F_{p,q},F_{r,s}\right]=-\delta_{q,r-1}F_{p,s}+\delta_{p,s+1}F_{r,q} \\
&\left[E_{p,q},F_{p,q}\right]=H_{p,q} & \\
& \left[E_{p,q},F_{p,s}\right]=-F_{q+1,s} \quad \mbox{if $s>q$} &\quad 
\left[E_{p,q},F_{p,s}\right]=
-E_{s+1,q}\quad \mbox{if  $q>s$}\\
& \left[E_{p,q},F_{r,q}\right]=E_{p,r-1} \quad \mbox{if $p<r$} &\quad  
\left[E_{p,q},F_{r,q}\right]=
F_{r,p-1}\quad \mbox{if  $p>r$}\\
&\left[H_i,P_j\right]=\delta_{ij}P_j-\delta_{i,j-1}P_j&\\
&\left[E_{p,q},P_j\right]=\delta_{q,j-1}P_{p}&\quad 
 \left[F_{p,q},P_j\right]=\delta_{p,j}P_{q+1}.\\
\end{array}
\end{equation}
\end{prop}
In \cite{DR1} Douglas and Repka have classified all the embeddings of the 
perfect but not simple Lie algebra $\mathfrak{sl}(n+1)\ltimes \mathbbm{C}^{n+1}$ 
in the simple Lie algebra $\mathfrak{sl}(n+2)$. 
Using for the  Lie algebra $\mathfrak{sl}(n+2)$ the same notion of above but 
with the indices running from $1$ to $n+1$ their results may summarized as 
\begin{theorem}\label{authomorphism}\cite{DR1}  There are,  up to inner 
automorphism, two inequivalent embeddings $\Phi$, $\Theta$ of \\
$\mathfrak{sl}(n+1)\ltimes \mathbbm{C}^{n+1}$ in the simple Lie algebra 
$\mathfrak{sl}(n+2)$:
\begin{equation}
\label{emmb1}
\begin{array}{ll}
&\Phi: \mathfrak{sl}(n+1)\ltimes \mathbbm{C}^{n+1} \xrightarrow{\makebox[2cm]{}}
 \mathfrak{sl}(n+2)\\
&\Phi(H_i)=A_{i+1,i+1}-A_{i+2,i+2}\mbox{ (the element $H_{i+1}$ of 
$\mathfrak{sl}(n+2)$)}\quad 1\leq i\leq  n\\
&\Phi(E_{p,q}) = A_{p+1,q+2}\mbox{ (the element $E_{p+1,q+1}$ of 
$\mathfrak{sl}(n+2)$)} 
 \quad 1\leq p\leq q\leq n\\
&\Phi(F_{p,q})= A_{q+2,p+1}\mbox{ (the element $F_{p+1,q+1}$ of 
$\mathfrak{sl}(n+2)$)} 
 \quad 1\leq p\leq q\leq n\\
&\Phi(P_{j})= A_{j+1,1}\mbox{ (the element $F_{1,j}$ of $\mathfrak{sl}(n+2)$)} 
 \quad 1\leq j\leq n+1\\
\end{array}
\end{equation} 
and 
\begin{equation}
\label{emmb2}
\begin{array}{ll}
&\Theta: \mathfrak{sl}(n+1)\ltimes \mathbbm{C}^{n+1} 
\xrightarrow{\makebox[2cm]{}} 
 \mathfrak{sl}(n+2)\\
&\Theta(H_i)=  -A_{i+1,i+1}+A_{i+2,i+2}\mbox{ (the element $-H_{i+1}$ of 
$\mathfrak{sl}(n+2)$)}\quad 1\leq i\leq  n\\
&\Theta(E_{p,q})=  -A_{q+2,p+1}\mbox{ (the element $-F_{q+1,p+1}$ of 
$\mathfrak{sl}(n+2)$)} 
 \quad 1\leq p\leq q\leq n\\
&\Theta(F_{p,q})= -A_{p+1,q+2}\mbox{ (the element $-E_{p+1,q+1}$ of 
$\mathfrak{sl}(n+2)$)} 
 \quad 1\leq p\leq q\leq n\\
&\Theta(P_{j}) = A_{1,j+1}\mbox{ (the element $E_{1,j}$ of 
$\mathfrak{sl}(n+2)$)} 
 \quad 1\leq j\leq n+1.\\
\end{array}
\end{equation} 
\end{theorem}
\begin{remark}\label{remb1} Both  embeddings  of $\mathfrak{sl}(n+1)\ltimes 
\mathbbm{C}^{n+1}$ into $\mathfrak{sl}(n+2)$,   $\Phi$, 
and $\Theta$ define on the  
space $\mathbbm{C}^{n+2}$ a  structure of $\mathfrak{sl}(n+1)\ltimes 
\mathbbm{C}^{n+1}$--module which we will denote respectively $\mathbbm{C}^{n+2}_\Phi$ and $\mathbbm{C}^{n+2}_\Theta$. Both modules are true and indecomposable.
The  Jordan--H\"older series of the module associated with the embedding 
$\Phi$ is 
\begin{equation}
\label{JH1}
0 = V_0 \subset V_1 \subset V_2=\mathbbm{C}^{n+2}_\Phi\qquad \mbox{with 
$\dim(V_2/V_1)=n+1$}.
\end{equation}
While that of the module associated to the embedding $\Theta$ is    
\begin{equation}
\label{JH2}
0 = W_0 \subset W_1 \subset W_2=\mathbbm{C}^{n+2}_\Theta\qquad \mbox{with 
$\dim(W_2/W_1)=1$.}
\end{equation}
\end{remark}
As a consequence of such remark we have 
\begin{prop}\label{propemb} Let $\Xi$ be the automorphims of 
the rooth space $\Delta$ of $\mathfrak{sl}(n+2)$ given by the relations 
$$
\Xi(\alpha_i)=\alpha_{n+2-i}\qquad 1\leq i\leq n+1.
$$ 
Denote also by  $\Xi$ the corresponding  automorphism of 
$\mathfrak{sl}(n+2)$ generated by the relations
$$
\Xi(H_{\alpha_i})=H_{\Xi(\alpha_i)}\quad \Xi(E_{\alpha_i})=E_{\Xi(\alpha_i)}
\quad \Xi(F_{\alpha_i})=F_{\Xi(\alpha_i)} \qquad 1\leq i\leq n =1
$$
on the Cartan basis of $\mathfrak{sl}(n+2)$. \par \noindent 
Then the embedding $\Xi\circ \Theta$, and $\Phi$ of 
$\mathfrak{sl}(n+1)\ltimes \mathbbm{C}^{n+1}$ in   $\mathfrak{sl}(n+2)$
are  equivalent. 
\end{prop}
{\bf Proof} It is easy to show that the map $\Xi$ is the restriction on 
$\mathfrak{sl}(n+2)$ of the automorphism of the space of the $(n+2)\times (n+2)$ 
complex matrices given by the equations:
$$
\Xi(A_{i,j})=(-1)^{j-i+1}A_{n+3-j,n+3-i}\qquad 1\leq i, j\leq n+2.
$$
This implies that 
$$
\begin{array}{ll}
&\hskip -0.5truecm \Xi\circ\Theta(H_i)=-A_{n+1-i,n+1-i}+ A_{n+2-i,n+2-i}
\quad 1\leq i\leq n\\
&\hskip -0.5truecm \Xi\circ\Theta(E_{p.q})=(-1)^{q-p}A_{n+2-p,n+1-q},\quad  
\Xi\circ\Theta(F_{p.q})= (-1)^{q-p}A_{n+1-q,n+2-p}\quad 1\leq p\leq q\leq n \\
&\Xi\circ\Theta(P_j)= (-1)^{j+1}A_{n+2-j,n+2} \quad 1\leq j\leq n+1.
\end{array}
$$
Hence the Jordan--H\"older series of the $n+2$ dimensional 
$\mathfrak{sl}(n+1)\ltimes \mathbbm{C}^{n+1}$--module associated with this 
embedding is 
 $$
0 = Z_0 \subset Z_1 \subset Z_2=\mathbbm{C}^{n+2}_{\Theta\circ\Xi}\qquad \mbox{with 
$\dim(Z_2/Z_1)=n+1$}.
$$
This latter equation together with equation (\ref{JH1}) show the equivalence 
of the embedding $\Xi\circ\Theta$ and $\Phi$.
\endpf 
 \section{Cyclic  $\mathfrak{sl}(n+1)\ltimes \mathbbm{C}^{n+1}$--modules}
Let $\mathfrak{g}$ be a Lie algebra, a module $M$ of $\mathfrak{g}$ is said  
indecomposable if it can not decomposed in the direct sum of two non trivial 
modules. An important class of indecomposable modules is given by 
\begin{defi}\label{cyclic} A module $M$ of a Lie algebra $\mathfrak{g}$ is said 
cyclic if  there exists an element $v\in M$ which generates $M$, i.e.,   
$M=\mathcal{U}(\mathfrak{g})v$.
\end{defi}
Obviously cyclic modules are indecomposable.  
Beyond the irreducible modules examples of cyclic modules are the string modules 
of the Euclidean Lie algebra in two dimensions \cite{RdG}, the modules of  the 
Diamond Lie algebra obtained by embedding it in $\mathfrak{sl}(3,\mathbbm{C})$ 
\cite{CMS}, the uniserial  modules \cite{CF1},  and the modules of $\mathfrak{sl}(2,\mathbbm{C})\ltimes \mathbbm{C}^2$ classified in \cite{Pi}.\par
Let  $\mathfrak{p}$ be   the subalgebra  of $\mathfrak{sl}(n+1)\ltimes
\mathbbm{C}^{n+1}$--modules, spanned by the elements $\{P_i\}$ $1\leq i\leq n+1$ 
(\ref{omega1basis}), then since $\mathfrak{p}$  is  an ideal in 
$\mathfrak{sl}(n+1)\ltimes \mathbbm{C}^{n+1}$ it holds 
 \begin{lem}\label{u(P)} Let $M$ be a  $\mathfrak{sl}(n+1)\ltimes 
\mathbbm{C}^{n+1}$--module   and $N$ be a 
$\mathfrak{sl}(n+1)$--module  contained 
in $M$. Then   the space   $\mathcal{U}(\mathfrak{p})N$ is a $\mathfrak{sl}(n+1)$--submodule.
\end{lem}
\begin{lem}\label{nil}
In any finite dimensional  $\mathfrak{sl}(n+1)\ltimes \mathbbm{C}^{n+1}$--module $M$, the elements $P_i$,  $1\leq i\leq n+1$ 
act as nilpotent operators and there exists a positive integer $m$ such that 
$M$ is annihilated by any monomial  of the type $P_{n+1}^{a_{n+1}}\cdots 
P_1^{a_1}$ with $a_{n+1}+\cdots + a_1\geq m$.
\end{lem}
{\bf Proof} Viewed as $\mathfrak{sl}(n+1)$--modules,   $M$ can  be 
 decomposed in the direct sum of its weight spaces $M=\bigoplus V_\mu$.
Now from the commutations relations 
$\left[H_i,P_j\right]=\delta_{ij}P_j-\delta_{i,j-1}P_j$, $1\leq i\leq n$,
$1\leq j\leq n+1$ (\ref{comrel}) it follows that the elements $P_i$  act 
on the weight space $V_\mu$ as
$$
P_i V_\mu \xrightarrow{\makebox[2cm]{}} V_{\mu+\omega_{i} -\omega_{i-1}}
$$
and therefore on any finite dimensional module they  are nilpotent operators.
Finally since the operators $P_i$ commute among themselves any monomial of the 
type  $P_{n+1}^{a_{n+1}}\cdots P_1^{a_1}$ is a nilpotent operator, and this implies the second statement of the Lemma.
\endpf
\begin{cor}\label{irri} A $\mathfrak{sl}(n+1)\ltimes \mathbbm{C}^{n+1}$--module
 $M$ is irreducible if and only if $\mathfrak{p}$ acts on $M$ trivially and 
$M$ is an irreducible $\mathfrak{sl}(n+1)$--module $M$.
\end{cor}
This latter result actually is a consequence of the general fact that in  a
perfect Lie algebra $\mathfrak{g}$,  its  solvable
radical $\mathfrak{r}$  coincides with its nilpotent radical 
$\left[\mathfrak{g}, \mathfrak{r}\right]$ \cite{CF1}. \par 
Any   finite dimensional $\mathfrak{sl}(n+1)\ltimes \mathbbm{C}^{n+1}$--module 
$M$ can be decomposed in a direct sum of irreducible finite dimensional 
 $\mathfrak{sl}(n+1)$--modules:
\begin{equation}
\label{dec1}
M=\bigoplus_{\lambda\in \Lambda} \pi(\lambda)V(\lambda)
\end{equation}
where $\pi(\lambda)$ is the    multiplicity of $V(\lambda)$ in $M$.\par 
In order to  classify the cyclic $\mathfrak{sl}(n+1)\ltimes 
\mathbbm{C}^{n+1}$--modules we need  to study how the  irreducible finite dimensional 
 $\mathfrak{sl}(n+1)$--modules are ``intertwined'' by the action of the ideal  $\mathfrak{p}$.\par 
\begin{defi}\label{mskdef}
Let $\mu=\sum_{i=1}^n \mu_i\omega_i \in P$ be a integral weight of
$\mathfrak{sl}(n+1)$, for any pair of integer numbers $(k,s)$ with $1\leq k\leq s\leq n$ 
we set 
\begin{equation}
\label{msk}
\mu_{k,s}=(-1)\sum_{i=k}^s(\mu_i+1)=-(\mu_s+\mu_{s-1}+\dots +\mu_k+s-k+1).
\end{equation}
\end{defi} 
By the definition of dominant weight and that of $\mu_{k,s}$ follows 
\begin{lem}\label{mk} If $\mu$ is a dominant weight then for any pair of positive integer  numbers $(k,s)$ with $1\leq k\leq s\leq n$ it holds  $\mu_{k,s}>0$.
\end{lem}
\begin{defi}\label{PmuFpq} For any integral   weight 
$\mu \in P$, and any integer $s$,  
$1\leq s\leq n$   we set    
 \begin{equation}
\label{pms}
P_i^{(\mu(s))}=\left(\prod_{k=1}^{i-1} \mu_{k,s-1}\right) P_i
\quad 1\leq i\leq s.
\end{equation} 
Observe that if $v_\mu$ is a weight vector of weight 
$\mu$ then $P_i^{(\mu(s))}(v_\mu)$ is, if different from zero, a weight vector of weight 
$\nu=\mu+\omega_{i}-\omega_{i-1}$, $1\leq i\leq n$.\par 
Further  let for any integral  weight $\mu \in P$, and any integer $s$, we set 
\begin{equation}
\label{fijms}
\begin{array}{lll}
 F_{i,j}^{(\mu(s))}&=\left(\displaystyle\prod_{k=i+1}^j \mu_{k,s}\right) 
F_{ij} &\quad 1\leq 
i\leq j\leq s\\
 F_{i,i}^{(\mu(s))}&= F_{i}^{(\mu(s))}=F_{i} &\quad 1\leq i\leq n.
\end{array}
\end{equation}
\end{defi}
Using this formalism let us give 
\begin{defi}\label{qr} For any  integral weight $\mu \in P$, any integer numbers  
$s$, $k$ with  $1\leq s\leq n$, $0\leq k< s$ let 
$Q^{(\mu(s))}_{s-k}$ be the element 
of  $\mathcal{U}(\mathfrak{sl}(n+1))$ defined by the recurrence relations:
\begin{equation}
\label{qmske}
\begin{array}{lll}
Q^{(\mu(s))}_{s+1}&=I & \\
Q^{(\mu(s))}_{s-k}&=\displaystyle\sum_{l=0}^k  
Q^{(\mu(s))}_{s-l+1}F^{(\mu(s))}_{s-k,s-l} 
&\quad 0\leq k\leq s-1
\end{array}
\end{equation}   
where $I$ is the identity of $\mathcal{U}(\mathfrak{sl}(n+1))$ 
\end{defi}
For example for $s=3$ we have: 
$$
\begin{array}{ll}
Q^{(\mu(3))}_{4}&=I  \\
Q^{(\mu(3))}_{3}&=F^{(\mu(3))}_{3}\\
Q^{(\mu(3))}_{2}&=F^{(\mu(3))}_{2,3}+ F^{(\mu(3))}_{3}F^{(\mu(3))}_{2}\\
Q^{(\mu(3))}_{1}&=F^{(\mu(3))}_{1,3}+ F^{(\mu(3))}_{3}F^{(\mu(3))}_{1,2}
+ F^{(\mu(3))}_{2,3}F^{(\mu(3))}_{1}+ 
F^{(\mu(3))}_{3}F^{(\mu(3))}_{2}F^{(\mu(3))}_{1}.
\end{array}
$$ 
\begin{defi}\label{rr} For any  integral weight $\mu \in P$, any  integer 
numbers  $s$, $i$,  $k$,  with $1\leq s\leq n$, $1\leq i\leq  s$,  $0\leq k \leq 
s-i+1  $,  let $R^{(\mu(s))}_{i,k}$ be the element 
of  $\mathcal{U}(\mathfrak{sl}(n+1))$ defined by the recurrence relations:
\begin{equation}
\label{rreq}
\begin{array}{llll}
R^{(\mu(s))}_{i,0}&=I &\quad 1\leq i\leq s  &\\
R^{(\mu(s))}_{i,k}&=\displaystyle\sum_{l=0}^{k-1}  F^{(\mu(s))}_{i+l,i+k-1}R^{(\mu(s))}_{i,l} 
&\quad 1\leq i\leq s &\quad 1\leq k\leq s-i+1
\end{array}
\end{equation} 
 where $I$ is the identity of $\mathcal{U}(\mathfrak{sl}(n+1))$. 
\end{defi}
For example for $s>3$ and $1\leq i\leq s$ we have:
$$
\begin{array}{ll}
R^{(\mu(s))}_{i,0}&=I \\
R^{(\mu(s))}_{i,1}&=F^{(\mu(s))}_{i}\\
R^{(\mu(s))}_{i,2}&=F^{(\mu(s))}_{i,i+1}+ F^{(\mu(s))}_{i+1}F^{(\mu(s))}_{i}\\
R^{(\mu(s))}_{i,3}&=F^{(\mu(s))}_{i,i+2}+ F^{(\mu(s))}_{i+1,i+2}F^{(\mu(3))}_{i}
+ F^{(\mu(s))}_{i+2}F^{(\mu(s))}_{i,i+1}+ 
F^{(\mu(s))}_{i+2}F^{(\mu(s))}_{i+1}F^{(\mu(s))}_{i}.
\end{array}
$$
Using the definitions of $P^{(\mu(s))}_i$,  $Q^{(\mu(s))}_k$  
$R^{(\mu(s))}_{i,k}$ it is easy to show
\begin{prop} \label{Qp} Let $\mu \in P $ be an integral weight. Let $s$, $k$ be   
integer numbers,  such that $1\leq s\leq n$, $0\leq k< s$. 
Then 
\begin{enumerate}
\item 
\begin{equation}
\label{qfij}
Q^{(\mu(s))}_k=\sum_{l=k}^s Q^{(\mu(s))}_{l+1} F^{(\mu(s))}_{k,l}. 
\end{equation}
\item $Q^{(\mu(s))}_k$ is a polynomial in $F^{(\mu(s))}_{p,q}$ with 
$k\leq p\leq q\leq s$, whose powers are  at most one. 
\item 
\begin{equation}
\label{2ks}
Q^{(\mu(s))}_k=\sum_{l=1}^{2^{s-k}} M^{(\mu(s))}_{k,l}\qquad 1\leq k\leq s
\end{equation}
with  $M^{(\mu(s))}_{k,l}$  monomials in 
$F^{(\mu(s))}_{p,q}$, $k\leq p\leq q\leq s$ such that for any pair of integer 
numbers $(i,j)$, with  
$k\leq i\leq j\leq s$ there exists at least one  monomial $M^{(\mu(s))}_{k,l}$ 
such that  $F_{i,j}$ divides $M^{(\mu(s))}_{k,l}$.
\item For any $1\leq s\leq n$, any  $1\leq k\leq s$, and any $1\leq l\leq 2^{s-k}$ there exists an 
integer  $h$ with $ k\leq h\leq s$ such that $F_{h,s}$ divides 
$M^{(\mu(s))}_{k,l}$
\item If both $F_{p,q}$ and $F_{i,j}$, with $k\leq p\leq q\leq s$ and    
$k\leq i\leq j\leq s$ divide $M^{(\mu(s))}_{k,l}$ then 
either $p\leq q <i\leq j$ or $i\leq j <p\leq q$.
\item For any monomial $M^{(\mu(s))}_{k,l}$ in equation (\ref{2ks}) there exist 
integer numbers\\  $n\geq j_1\geq j_2\geq \cdots \geq j_h\geq 1$ such that
$$
M^{(\mu(s))}_{k,l}=\prod_{i=1}^{h+1}\left(\prod_{l=j_i+1}^{j_{i-1}}\mu_{l,s}\right)F_{j_1,s}F_{j_2,j_1-1}\cdots F_{j_h,j_{h-1}-1}F_{k,j_h-1}
$$
where  we have set as before $F_{j,0}=1$ and $j_0=s$.
\end{enumerate}
\end{prop}
{\bf Proof} 
\begin{enumerate}
\item It is nothing else but equation (\ref{qmske}) written in a more convenient 
way.
\item It is an immediate consequence of the very form of equation
 (\ref{qfij}) and of the definition of the elements $Q^{(\mu(s))}_k$, because the
elements $F^{(\mu(s))}_{p,q}$ with $p\leq q\leq s$ are not factors of the 
monomials which appear in the expression of $Q^{(\mu(s))}_k$ with $k>q$.
\item It follows by induction. For $k=s$ the claim is clear. Suppose now 
that the claim is true for 
$Q^{(\mu(s))}_h$ with $h>k$ then from   (\ref{qfij}) follows that it is true for 
$Q^{(\mu(s))}_k$,  because if the index $i$, in  $F^{(\mu(s))}_{i,j}$  is $i> k$ 
then it divides by induction at least  one of the monomials of $Q^{(\mu(s))}_h$ with $h>k$, while   if 
$i=k$,  $F^{(\mu(s))}_{k,j}$ appears explicitly in  (\ref{qfij}). 
\item It follows immediately by induction.  
\item Again it follows by induction. For $k=s$ it is clear. Now if the claim is 
true for $Q^{(\mu(s))}_k$ with $k>j$ then from (\ref{qfij}) and (2) follows that 
$F_{j,l}$ for any $l$ with  $j\leq l\leq s$ appears as a factor only in 
monomials  of $Q^{(\mu(s))}_j$ in which any other factor  $F_{p,q}$ has $j\leq l<p\leq q$.
\item It follows almost directly from the previous points.
\end{enumerate}
\endpf
\begin{prop} \label{Rp} Let $\mu \in P $ an integral weight. Let $s$,  $i$,
 $h$  integer numbers, such that  $1\leq s\leq n$,  $1\leq i\leq  s$, 
$1\leq h\leq  s-i+1$. Then 
\begin{enumerate}
\item 
\begin{equation}
\label{rfij}
R^{(\mu(s))}_{i,h}=\sum_{l=0}^{h-1} F^{(\mu(s))}_{i+l,i+h-1} R^{(\mu(s))}_{i,l}.
\end{equation}
\item 
$R^{(\mu(s))}_{i,h}$ is a polynomial in $F^{(\mu(s))}_{p,q}$ with 
$i\leq p\leq q\leq i+h-1$, whose power is at most one. 
\item 
$$
R^{(\mu(s))}_{i,h}=\sum_{l=1}^{2^{h-1}}N^{(\mu(s))}_{i,h,l}
$$
 with   $N^{\mu(s))}_{i,h,l}$  monomials in 
$F^{(\mu(s))}_{p,q}$, $i\leq p\leq q\leq i+h-1$ such that for any pair $(j,m)$, 
$i\leq j\leq m\leq i+h-1$ there exist at least one  monomial
$N^{(\mu(s))}_{i,h,l}$  such that  $F_{j,m}$ divides $N^{(\mu(s))}_{i,h,l}$.
\item For any $1\leq s\leq n$,  any $1\leq i\leq s$, any   $1\leq h\leq s-i+1$, and 
any $1\leq l\leq 2^{h-1}$ there exists an 
integer  $m$ with $ k\leq m\leq i+h-1$ such that $F_{m,i+h-1}$ divides 
$N^{(\mu(s))}_{i,h,l}$
\item If both $F_{p,q}$ and $F_{j,m}$, $k\leq p\leq q\leq i+h-1$,   
$k\leq j\leq m\leq i+h-1$ divide $N^{(\mu(s))}_{i,h,l}$ then 
either $p\leq q <j\leq m$ or $j\leq m<p\leq q$.
\item $R^{(\mu(s))}_{i,s-i+1}=Q^{(\mu(s))}_i,\qquad 1\leq i\leq s$.
\item If $v_\mu$ is a weight vector of weight $\mu$ then the vector  
$R^{(\mu(s-1))}_{i,h}P^{(\mu(s))}_i(v_\mu)$, if not equal to zero,  is a weight vector of weight
$\mu +\omega_{i+h}-\omega_{i+h-1}$. 
\end{enumerate}
\end{prop}
{\bf Proof} Mutata mutandis, for the first five  points we can argue as in the 
previous Proposition. Point 6.   follows immediately from the Definitions \ref{Qp} and \ref{Rp}; point 7. from the  
properties  of the operators $P_i$ and the commutation rules (\ref{comrel}).
\endpf
Furthermore using the commutation rules of 
$\mathfrak{sl}(n+1)\ltimes \mathbbm{C}^{n+1}$,  the Definitions 
\ref{comrel}, the equations (\ref{pms}) and (\ref{fijms}) and the above
Propositions \ref{qr} \ref{rr},  we have 
\begin{prop} \label{EFP} Let $\mu \in P $ be  an integral weight, $1\leq s\leq n$. Then for any  $1\leq j\leq s$ we have
\begin{equation}
\label{EFPeq}
\begin{array}{llll}
\left[ E_j,F^{(\mu(s))}_{j,l}\right]&=-\mu_{j+1,s}F^{(\mu(s))}_{j+1,l} \ \ 
1\leq j< l
&\quad   \left[ E_j,F^{(\mu(s))}_{l,j}\right]=\mu_{j,s}
F^{(\mu(s))}_{l,j-1}  &\ \ \ \
1\leq l< j\\
\left[E_j,F^{(\mu(s))}_j\right]&=H_j   &\quad 
 \left[ E_j,Q^{(\mu(s))}_{l}\right]=0    &\ \     j< l,\ j>s\\
\left[ E_j,R^{(\mu(s))}_{i,h}\right]&=0 \ \  i>j,\ i+h\leq j  &\quad    
\left[ E_j,P^{(\mu(s))}_{i}\right]
=\delta_{j,i-1}\mu_{j,s-1}P^{(\mu(s))}_{j} 
&\ \   1\leq i\leq s. 
\end{array}
\end{equation}
\end{prop}
Armed with these facts we can  prove the 
\begin{theorem}\label{hwv} Let $V$ be a  
$\mathfrak{sl}(n+1)\ltimes \mathbbm{C}^{n+1}$--module, and let $v_\mu\in V$ be a 
highest weight vector of  
$\mathfrak{sl}(n+1)$ of weight $\mu$. Then for 
$1\leq i\leq n+1$ the element  
\begin{equation}
\label{hwveq}
\varphi_{i}(v_\mu)=P^{(\mu(i))}_{i}(v_\mu)+\sum_{k=1}^{i-1}Q^{(\mu(i-1))}_kP^{(\mu(i))}_{k}(v_\mu)\qquad 
\end{equation}
is  either zero or   a highest weight vector of   $\mathfrak{sl}(n+1)$ of 
weight  $\mu+\omega_{i}-\omega_{i-1}$. 
\end{theorem}
{\bf Proof} It is enough to show that for any $j$ with $1\leq j\leq n$ 
\begin{equation}
\label{ejphi}
E_j\ \varphi_{i}(v_\mu)=0\qquad 1\leq i\leq n+1.
\end{equation}
Now from the Propositions  \ref{Qp}, \ref{Rp} and \ref{EFP} and the definition of $\varphi_{i}(v_\mu)$  it is obvious that
$$
E_j\  \varphi_{i}(v_\mu)=0 \qquad \mbox{for } j\geq i.
$$ 
Hence we need only to consider the elements $E_j\varphi_{i}(v_\mu)$ 
with $j\leq i-1$.
Using (\ref{hwveq}) we can write $\varphi_i(v_\mu)$ in this case as
$$
\varphi_i(v_\mu)=\sum_{k=j+2}^i Q^{(\mu(i-1))}_k P^{(\mu(i))}_k(v_\mu)+
\sum_{k=1}^{j+1}Q^{(\mu(i-1))}_k P^{(\mu(i))}_k(v_\mu).
$$
Again using  propositions  \ref{Qp}, \ref{Rp} and \ref{EFP}  we have
$$
E_j\left(\sum_{k=j+2}^i Q^{(\mu(i-1))}_k P^{(\mu(i))}_k(v_\mu)\right)=0.
$$
So we need only to consider the action of $E_j$ on 
$\sum_{k=1}^{j+1}Q^{(\mu(i-1))}_k P^{(\mu(i))}_k(v_\mu)$.\par 
From the formulas (\ref{EFPeq}) it follows that if $E_j$ annihilates 
$\sum_{k=1}^{j+1}Q^{(\mu(i-1))}_k P^{(\mu(i))}_k(v_\mu)$ then it must 
hold
\begin{align}
&E_j \left(Q^{(\mu(i-1))}_{j+1} P^{(\mu(i))}_{j+1}(v_\mu)+ Q^{(\mu(i-1))}_{j} 
P^{(\mu(i))}_{j}(v_\mu)\right)=0 \label{ejh1}\\
\nonumber \\
&E_j \left(Q^{(\mu(i-1))}_k P^{(\mu(i))}_{k}(v_\mu)\right)=0  \quad 1\leq k<
j.\label{ejh2}
\end{align} 
Let us start by considering equation (\ref{ejh1}). We have two possible cases:
either $i=j+1$ or $i>j+1$. In first case  (\ref{ejh1}) becomes
$$
E_j \left(P^{(\mu(j+1))}_{j+1}(v_\mu)+ F^{(\mu(j))}_{j} 
P^{(\mu(j+1))}_{j}(v_\mu)\right)=0
$$
and we have using equations (\ref{EFPeq})
$$
\begin{array}{ll}
&E_j \left(P^{(\mu(j+1))}_{j+1}(v_\mu)+ F^{(\mu(j))}_{j} 
P^{(\mu(j+1))}_{j}(v_\mu)\right)\\
&=E_j\left(P^{(\mu(j+1))}_{j+1}(v_\mu)\right)+
\left[E_j,F^{(\mu(j))}_{j}\right]P^{(\mu(j+1))}_{j}(v_\mu)
=-(\mu_j+1) P^{(\mu(j+1))}_{j}(v_\mu)+H_jP^{(\mu(j+1))}_{j}(v_\mu)\\
&=-(\mu_j+1) P^{(\mu(j+1)}_{j}(v_\mu)+(\mu_j+1)P^{(\mu(j+1))}_{j}(v_\mu)=0.
\end{array}
$$
While if $i>j+1$, 
equation (\ref{ejh1}) becomes
$$
E_j \left(Q^{(\mu(i-1))}_{j+1}P^{(\mu(i))}_{j+1}(v_\mu)+ Q^{(\mu(i-1))}_{j} 
P^{(\mu(i))}_{j}(v_\mu)\right)=0.
$$
Now using (\ref{qfij}) we have 
$$
\begin{array}{ll} 
Q^{(\mu(i-1))}_{j+1}&=\sum_{l=j+1}^{i-1}Q^{(\mu(i-1))}_{l+1}F^{(\mu(i-1))}_{j+1,l}\\
Q^{(\mu(i-1))}_{j}&=\sum_{l=j}^{i-1}Q^{(\mu(i-1))}_{l+1}F^{(\mu(i-1))}_{j,l}\\
&=\sum_{l=j+1}^{i-1}Q^{(\mu(i-1))}_{l+1}F^{(\mu(i-1))}_{j,l}+
\sum_{l=j+1}^{i-1}Q^{(\mu(i-1))}_{l+1}F^{(\mu(i-1))}_{j+1,l}F^{(\mu(i-1))}_j
\end{array}
$$
therefore using (\ref{msk}), (\ref{EFPeq}) and Proposition \ref{Rp}
$$
\begin{array}{ll} 
&E_j \left(Q^{(\mu(i-1))}_{j+1}P^{(\mu(i))}_{j+1}(u_\mu)+ Q^{(\mu(i-1))}_{j} 
P^{(\mu(i))}_{j}(v_\mu)\right)\\
&=\sum_{l=j+1}^{i-1}Q^{(\mu(i-1))}_{l+1}F^{(\mu(i-1))}_{j+1,l}
\left[E_j,P^{(\mu(i))}_{j+1}\right](v_\mu)+
\sum_{l=j+1}^{i-1}Q^{(\mu(i-1))}_{l+1}\left[E_j,F^{(\mu(i-1))}_{j,l}\right]
P^{(\mu(i))}_{j}(v_\mu)\\
&+\sum_{l=j+1}^{i-1}Q^{(\mu(i-1))}_{l+1}F^{(\mu(i-1))}_{j+1,l}
\left[E_j,F^{(\mu(i-1))}_j\right]P^{(\mu(i))}_{j}(v_\mu)\\
&=\mu_{j,i-1}\sum_{l=j+1}^{i-1}Q^{(\mu(i-1))}_{l+1}F^{(\mu(i-1))}_{j+1,l}
P^{(\mu(i))}_{j}(v_\mu)
-\mu_{j+1,i-1}\sum_{l=j+1}^{i-1}Q^{(\mu(i-1))}_{l+1}F^{(\mu(i-1))}_{j+1,l}
P^{(\mu(i))}_{j}(v_\mu)\\
&+(\mu_j+1)\sum_{l=j+1}^{i-1}Q^{(\mu(i-1))}_{l+1}F^{(\mu(i-1))}_{j+1,l}
P^{(\mu(i))}_{j}(v_\mu)=0.
\end{array}
$$
It remains to prove  equation (\ref{ejh2}). This requires some more 
formalism.\par 
Recall that from equation (\ref{2ks})  $Q_k^{(\mu(s))}$ for $1\leq k \leq s$ is the sum of monomials $M^{(\mu(s))}_{k,l}$ with $1\leq l\leq 2^{s-k}$.   
For any $j$, $k\leq j\leq s$   let $J^+_{k,s}(j)$ be the set:
$$
J^+_{k,s}(j)=\left\{ 1\leq l\leq 2^{s-k}\bigm| \  \exists\ h \  j\leq h\leq s \
\mbox{such that } F_{j,h} \mbox{ divides }  M^{(\mu(s))}_{k,l} \  \right\}
$$
i.e., the set of all $l$, $1\leq l\leq 2^{s-k}$ such that the monomial
$M^{(\mu(s))}_{k,l}$ contains  a factor of type $F_{j,h}$ with $j\leq h\leq s$.
Let $Q^{(\mu(s))^+}_{k,j}$ be  defined by
\begin{equation}
\label{qkjp}
Q^{(\mu(s))^+}_{k,j}=\sum_{l\in J^+_{k,s}(j)}  M^{(\mu(s))}_{k,l},
\end{equation}
 then from proposition \ref{Qp} anf \ref{Rp} follows almost immediately that
for any  $1\leq j\leq n$
\begin{equation}
\label{qkj}
Q^{(\mu(s))^+}_{k,j}=Q^{(\mu(s))}_jR^{(\mu(s))}_{k,j-k}.
\end{equation} 
Now from (\ref{qfij}) and (\ref{rfij}) we have
$$
\begin{array}{ll}
Q^{(\mu(s))}_j&=\sum_{l=j}^s  Q^{(\mu(s))}_{l+1}F^{(\mu(s))}_{j,l}\\
R^{(\mu(s))}_{k,j-k}&=\sum_{l=0}^{j-k-1}  F^{(\mu(s))}_{k+l,j-1} 
R^{(\mu(s))}_{k,l}\\
\end{array}
$$
using these latter equations in (\ref{qkj}) we have
\begin{equation}
\label{qkj1}
Q^{(\mu(s))^+}_{k,j}=\sum_{l=j}^s\sum_{h=0}^{j-k-1}
Q^{(\mu(s))}_{l+1}F^{(\mu(s))}_{j,l}F^{(\mu(s))}_{k+h,j-1}R^{(\mu(s))}_{k,h}.
\end{equation}
We still have to take into account those monomials in the expression of
$Q^{(\mu(s))}_k$ which are divisible for $F_{h,j}$ with $k\leq  h < j$
(those  with $h=j$ has been already considered).\par
As above, for any $j$, $k\leq j\leq s$   let $J^-_{s,k}(j)$ be the set:
$$
J^-_{k,s}(j)=\left\{ 1\leq l\leq 2^{s-k}\bigm| \  \exists\ h, 
\ k\leq h <j\mbox{ such that $F_{h,j}$  divides   $M^{(\mu(s))}_{k,l}$}  \right\}
$$
i.e., the set of all $l$, $1\leq l\leq 2^{s-k}$ such that the monomial
$M^{(\mu(s))}_{k,l}$ contains  a factor of type $F_{h,j}$ with $k\leq h<j$.
Let $Q^{(\mu(s))^-}_{k,j}$ be  defined by
\begin{equation}
\label{qkjm}
Q^{(\mu(s))^-}_{k,j}=\sum_{l\in J^-_{k,s}(j)}  M^{(\mu(s))}_{k,l}.
\end{equation}
Arguing as in the previous case  we have that 
\begin{equation}
\label{qrkj}
Q^{(\mu(s))^-}_{k,j}=Q^{(\mu(s))}_{j+1}\left(R^{(\mu(s))}_{k,j-k+1}-
\mbox{ the monomials in $R^{(\mu(s))}_{k,j-k+1}$ which are  divisible by 
$F_j$}\right).
\end{equation}
Now from Proposition \ref{Rp} and equation  (\ref{rfij}) it follows that 
$$
\left(R^{(\mu(s))}_{k,j-k+1}-
\mbox{ the monomials in $R^{(\mu(s))}_{k,j-k+1}$ which are  divisible by 
$F_j$}\right)=\sum_{h=0}^{j-k-1} F_{k+h,j}^{(\mu(s))} R^{(\mu(s))}_{k,h}.
$$
Therefore using this latter equation together with (\ref{qfij}) in (\ref{qrkj}) 
we have  
\begin{equation}
\label{qrkj1}
Q^{(\mu(s))^-}_{k,j}=
\sum_{l=j+1}^{s}\sum_{h=0}^{j-k-1}Q^{(\mu(s))}_{l+1}F^{(\mu(s))}_{j+1,l}
F^{(\mu(s))}_{k+h,j}R^{(\mu(s))}_{k,h}.
\end{equation}
Now, for what said above, since obviously $\left[ E_j, 
M^{(\mu(i-1))}_{k.l}\right]=0$ if $M^{(\mu(i-1))}_{k.l}$ is  not divisible by 
$F_{h,j}$, $k\leq h\leq j$ or $F_{j,h}$ $j\leq h\leq i-1$ we have only to show 
that
$$
E_j(Q^{(\mu(i-1))^+}_{k,j}P^{(\mu(i))}_k(v_\mu)
+Q^{(\mu(i-1))^-}_{k,j}P^{(\mu(i))}_k(v_\mu))=0.
$$
Taking into account equations (\ref{qkj1}),  (\ref{qrkj1}) and Proposition 
\ref{Rp} we have 
$$
\begin{array}{ll}
&E_j\left(\sum_{l=j+1}^{i-1}\sum_{h=0}^{j-k-1} Q^{(\mu(i-1))}_{l+1}
F^{(\mu(i-1))}_{j,l}F^{(\mu(i-1))}_{k+h,j-1}R^{(\mu(i-1))}_{k,h}P^{(\mu(i))}_k(v_\mu)\right.\\
&+\left.\sum_{h=0}^{j-k-1} Q^{(\mu(i-1))}_{j+1}
F^{(\mu(i-1))}_jF^{(\mu(i-1))}_{k+h,j-1}
R^{(\mu(i-1))}_{k,h}P^{(\mu(i))}_k(v_\mu)\right.\\
&\left.+\sum_{l=j+1}^{i-1}\sum_{h=0}^{j-k-1} Q^{(\mu(i-1))}_{l+1}
F^{(\mu(i-1))}_{j+1,l}F^{(\mu(i-1))}_{k+h,j}R^{(\mu(i-1))}_{k,h}P^{(\mu(i))}_k(v_\mu)\right)\\
&=E_j\left(\sum_{l=j+1}^{i-1}\sum_{h=0}^{j-k-1} Q^{(\mu(i-1))}_{l+1}
F^{(\mu(i-1))}_{j,l}F^{(\mu(i-1))}_{k+h,j-1}R^{(\mu(i-1))}_{k,h}P^{(\mu(i))}_k(v_\mu)\right.\\
&\left.+\sum_{l=j+1}^{i-1}\sum_{h=0}^{j-k-1} Q^{(\mu(i-1))}_{l+1}
F^{(\mu(i-1))}_{j+1,l}F^{(\mu(i-1))}_jF^{(\mu(i-1))}_{k+h,j-1}
R^{(\mu(i-1))}_{k,h}P^{(\mu(i))}_k(v_\mu)\right.\\
&\left.+\sum_{l=j+1}^{i-1}\sum_{h=0}^{j-k-1} Q^{(\mu(i-1))}_{l+1}
F^{(\mu(i-1))}_{j+1,l}F^{(\mu(i-1))}_{k+h,j}R^{(\mu(i-1))}_{k,h}P^{(\mu(i))}_k(v_\mu)\right)\\
&=\sum_{l=j+1}^{i-1}\sum_{h=0}^{j-k-1} Q^{(\mu(i-1))}_{l+1}
\left[E_j,F^{(\mu(i-1))}_{j,l}\right]F^{(\mu(i-1))}_{k+h,j-1}
R^{(\mu(i-1))}_{k,h}P^{(\mu(i))}_k(v_\mu)\\
&+\sum_{l=j+1}^{i-1}\sum_{h=0}^{j-k-1} Q^{(\mu(i-1))}_{l+1}
F^{(\mu(i-1))}_{j+1,l}\left[E_j,F^{(\mu(i-1))}_j\right]F^{(\mu(i-1))}_{k+h,j-1}
R^{(\mu(i-1))}_{k,h}P^{(\mu(i))}_k(v_\mu)\\
&+\sum_{l=j+1}^{i-1}\sum_{h=0}^{j-k-1} Q^{(\mu(i-1))}_{l+1}
F^{(\mu(i-1))}_{j+1,l}\left[E_j,F^{(\mu(i-1))}_{k+h,j}\right]
R^{(\mu(i-1))}_{h,k}P^{(\mu(i))}_k(v_\mu)\\
&=-\mu_{j+1,i-1}\sum_{l=j+1}^{i-1}\sum_{h=0}^{j-k-1} Q^{(\mu(i-1))}_{l+1}
F^{(\mu(i-1))}_{j+1,l}F^{(\mu(i-1))}_{k+h,j-1}
R^{(\mu(i-1))}_{k,h}P^{(\mu(i))}_k(v_\mu)\\
&+(\mu_j+1)\sum_{l=j+1}^{i-1}\sum_{h=0}^{j-k-1} Q^{(\mu(i-1))}_{l+1}
F^{(\mu(i-1))}_{j+1,l}F^{(\mu(i-1))}_{k+h,j-1}R^{(\mu(i-1))}_{k,h}P^{(\mu(i))}_k(v_\mu)\\
&+\mu_{j,i-1}\sum_{l=j+1}^{i-1}\sum_{h=0}^{j-k-1} Q^{(\mu(i-1))}_{l+1}
F^{(\mu(i-1))}_{j+1,l}F^{(\mu(i-1))}_{k+h,j-1}R^{(\mu(i-1))}_{k,h}P^{(\mu(i))}_k(v_\mu)=0.
\end{array}
$$
\endpf
\begin{remark}
The first elements $\varphi_i(v_\mu)$ are
$$
\begin{array}{ll}
\varphi_1(v_\mu)&=P_1 v_\mu\\
\varphi_2(v_\mu)&=-(\mu_1+1)P_2 v_\mu+F_1P_1v_\mu\\
\varphi_3(v_\mu)&=(\mu_2+1)(\mu_2+\mu_1+2)P_3 v_\mu-(\mu_2+\mu_1+2)F_2P_2v_\mu\\
&-(\mu_2+1)F_{1,2}P_1v_\mu +F_2F_1P_1v_\mu\\
\varphi_4(v_\mu)&=-(\mu_3+1)(\mu_3+\mu_2+2)(\mu_3+\mu_2+\mu_1+3)P_4 v_\mu\\
&+(\mu_3+\mu_2+2)(\mu_3+\mu_2+\mu_1+3)F_3P_3 v_\mu
+(\mu_3+1)(\mu_3+\mu_2+\mu_1+3)F_{2,3}P_2v_\mu\\
&-(\mu_3+\mu_2+\mu_1+3)F_3F_{2}P_2v_\mu 
+(\mu_3+1)(\mu_3+\mu_2+2)F_{1,3}P_1v_\mu\\
&-(\mu_3+\mu_2+2)F_{3}F_{1,2}P_1v_\mu-(\mu_3+1)F_{2,3}F_{1}P_1v_\mu+
F_3F_2F_1P_1v_\mu. 
\end{array}
$$
\end{remark}
\begin{cor}\label{corrphii} Let $V$ be a $\mathfrak{sl}(n+1)\ltimes 
\mathbbm{C}^{n+1}$--module,  $V(\mu)\subset V$ be an irreducible 
$\mathfrak{sl}(n+1)$--module of highest weight $\mu$  and highest vector $v_\mu$, then for any $1\leq i\leq n+1$ there exists a positive integer  
$k_i$ such that
\begin{enumerate}
\item $\varphi_i^{k_i}(v_\mu)=0$
\item For any $i$,  $1\leq i\leq n+1$ and any $k$,  $1\leq k< k_i$, 
the space $\mathcal{U}(\mathfrak{sl}(n+1))\varphi^k_i(v_\mu)$ is a finite 
irreducible $\mathfrak{sl}(n+1)$--module with  highest weight 
$\mu+k(\omega_i-\omega_{i-1})$ and highest weight vector $\varphi^k_i(v_\mu)$.
\end{enumerate}
\end{cor}
{\bf Proof} That any $\varphi_i^k(v_\mu)\neq 0$, $1\leq i\leq n+1$, $k\in 
\mathbbm{N}$, is a highest weight vector for $\mathfrak{sl}(n+1)$ of weight 
$\mu+k(\omega_i-\omega_{i-1})$ has been 
already proved in the previous Theorem \ref{hwv}.\par  
The existence of  a such $k_i$ for any $1\leq i\leq n+1$ follows form 
the fact that if different from zero $\varphi^k_i(v_\mu)$ and $\varphi^h_i(v_\mu)$ 
have different weight if $k\neq h$ and that the space $V$ is finite dimensional. 
The implication $\varphi^k_i(v_\mu)=0\Longrightarrow \varphi_i^h(v_\mu)=0$ 
if $k\leq h$, $1\leq i\leq n+1$ proves the second statement of the Corollary.
\endpf
The formulas (\ref{hwveq}) may be inverted, let us give  indeed:
\begin{defi}\label{Fphihat} For any integral dominant weight $\mu\in \Lambda$ any positive integer numbers $j$, $k$, $i$  with $1\leq j\leq k\leq n$ and  $1\leq i\leq k$ define
\begin{equation}\label{fijhateq}
\widehat{F}_{j.k}^{\ (\mu,i)}=\frac{1}{\vert \mu_{i,k}\ \vert}F_{j,k}
\end{equation}
where $\vert x\ \vert$ is the absolute value of $x$.
\end{defi}
\begin{defi}\label{Qhat} For any integral dominant weight $\mu\in \Lambda$, any 
positive integer numbers $s$, $i$, $k$ with $1\leq i\leq s\leq n$ and $i\leq 
k\leq s+1$ let $\widehat{Q}_k^{(\mu(s),i)}$ be the element of 
$\mathcal{U}(\mathfrak{sl}(n+1))$ defined by the recurrence relation:
 \begin{equation}\label{Qhateq}
\begin{array}{ll}
\widehat{Q}^{\ (\mu(s),i)}_{s+1}&=I\\
\widehat{Q}^{\ (\mu(s),i)}_{k}&=\displaystyle\sum_{l=k+1}^{s+1}
\widehat{Q}^{\ (\mu(s),i)}_{l}
\widehat{F}_{k,l-1}^{\ (\mu,i)}\qquad i\leq k<s+1.
\end{array}
\end{equation}
\end{defi}
Observe  that from Lemma \ref{mk} it follows that the elements 
$ \widehat{F}_{j.k}^{(\mu,i)}$ (\ref{fijhateq}) and the following elements 
$\widehat{\varphi}_i(v_\mu)$ (\ref{varphihat}) as well are well defined. 
Then it holds 
\begin{prop}\label{invphiP} Let $V$ be a  
$\mathfrak{sl}(n+1)\ltimes \mathbbm{C}^{n+1}$--module, and let $v_\mu\in V$ be a 
highest weight vector of  $\mathfrak{sl}(n+1)$ of integral dominant weight $\mu$. Let further $\widehat{\varphi}_i(v_\mu)$, $1\leq i\leq n+1$ be defined by
\begin{equation}
\label{varphihat}
\widehat{\varphi}_i(v_\mu)=\frac{1}{\prod_{k=1}^{i-1}\vert \mu_{k,i-1}\ \vert } 
\varphi_i(v_\mu)\qquad 1\leq i\leq n+1.\ 
\end{equation}
Then the elements $P_iv_\mu$, $1\leq i\leq n+1$ can be written as 
\begin{equation}
\label{invvarphi}
P_iv_\mu=(-1)^{i+1}\widehat{\varphi}_i(v_\mu)+\sum_{k=1}^{i-1} (-1)^{k+1}
\widehat{Q}^{\ (\mu(i-1),k)}_k\widehat{\varphi}_k(v_\mu)
\qquad 1\leq i\leq n+1.
\end{equation}
\end{prop}
 {\bf Proof} The transformation between the vectors $P_iv_\mu$ and the vectors
$\varphi_i(v_\mu)$ $i=1,\dots n+1$ given by the equations (\ref{hwveq}) can be 
viewed as an upper triangular matrix with entries on the main diagonal different 
from zero. Therefore it can be inverted.\par 
It remains to  show that the explicit form of such inverse is given by equations 
(\ref{invvarphi}). We shall argue by induction over the rank  of the simple Lie algebra $\mathfrak{sl}(n+1)$ (i.e., over the integer number $n$). 
For $n=1$ the claim  is obvious.  We shall show that  
the equations (\ref{invvarphi}) hold  for $\mathfrak{sl}(n+1)\ltimes  
\mathbbm{C}^{n+1}$ if they  hold  for $\mathfrak{sl}(n)\ltimes 
\mathbbm{C}^{n}$.\par 
A moment thought reveals that for $1\leq i\leq n$   the explicit form of 
the element $\varphi_i(v_\mu)$  in a  $\mathfrak{sl}(n)\ltimes  
\mathbbm{C}^{n}$--module  is the same of that  of 
the corresponding element in a 
$\mathfrak{sl}(n+1)\ltimes \mathbbm{C}^{n+1}$--modules.
Hence we need only to prove that equations (\ref{invvarphi}) hold when $i=n+1$. From the very same 
reason we also need only to show that 
\begin{equation}
\label{pn+1}
P_{n+1}v_\mu=\sum_{l=2}^{n+1} S_l(F_{p,q})\varphi_l(v_\mu)+
\widehat{Q}^{\ (\mu(n),1)}_1\widehat{\varphi}_1(v_\mu)
\qquad S_l(F_{p,q})\in \mathcal{U}(\mathfrak{n}^-)\quad 2\leq l\leq n+1.
\end{equation}
Using equation (\ref{hwveq}) with $i=n+1$ together (by induction Hypothesis) 
with  equations (\ref{invvarphi}) with $i\leq n$ and the point 6. of Proposition \ref{Qp} we obtain 
$$
\begin{array}{ll}
&\hskip-0.5cm P_{n+1}v_\mu=\sum_{l=2}^{n+1} S_l(F_{p,q})\varphi_l(v_\mu)\\
&+\hskip-0.8cm\displaystyle\sum_{n\geq k_1\geq k_2\geq \cdots \geq k_h\geq 1} 
\frac{1}{\displaystyle\prod_{l=1}^{n}\vert 
\mu_{l,n}\vert} \sum_{i=1}^{h+1} 
\left(\frac{(-1)^{i+1}\prod_{l=1}^{n}\vert\mu_{l,n}\vert} 
{\prod_{l=1}^{i}\vert \mu_{k_l,n}\vert \prod_{j=i}^{h}
\vert \mu_{1,k_{j-1}}\vert}
\right) 
F_{n,k_1}
F_{k_2,k_1-1}
\cdots F_{k_{h},k_{h-1}-1}F_{1,k_h-1}
\varphi_1(v_\mu)
\end{array}
$$ 
where the $S_l(F_{p,q})$,  $2\leq l\leq n+1$  are polynomials  in the variable $F_{p,q}$ with  $2\leq p\leq q\leq n$.\par\noindent  
It is easy to compute that 
$$
\frac{1}{\displaystyle\prod_{l=1}^{n}\vert 
\mu_{l,n}\vert} \sum_{i=1}^{h+1}
\left(\frac{(-1)^{i+1}\prod_{l=1}^{n}\vert\mu_{l,n} \vert}
{\prod_{l=1}^{i}\vert\mu_{k_l,n}\vert\prod_{j=i}^{h}
\vert\mu_{1,k_{j-1}}\vert}
\right)=
\frac{1}{\displaystyle\prod_{l=0}^h
\vert\mu_{1,k_l-1}\vert},\quad k_0=n+1.
$$
Now Definition of $\widehat{Q}^{\ \mu(n),1}_1$ \ref{Qhat} 
and Proposition \ref{Qp} yield 
$$
\widehat{Q}^{\ \mu(n),1}_1=\displaystyle\sum_{n\geq k_1\geq k_2\geq \cdots \geq k_h\geq 1} 
\frac{1}{\displaystyle\prod_{l=0}^h\vert \mu_{1,k_l-1}\vert}
F_{n,k_1}F_{k_2,k_1-1}\cdots F_{k_h,k_{h-1}-1}F_{1,k_h-1}
$$
which implies equation (\ref{pn+1}) as wanted.
\endpf 
\begin{remark}
The first examples of equations (\ref{invvarphi})  are
$$
\begin{array}{ll}
P_1(v_\mu)&=\varphi_1 (v_\mu)\\
P_2v_\mu&=-\frac{1}{\mu_1+1}\varphi_2(v_\mu)
+\frac{1}{\mu_1+1}F_1\varphi_1(v_\mu)\\
P_3(v_\mu)&=\frac{1}{(\mu_2+1)(\mu_2+\mu_1+2)}
\varphi_3 (v_\mu)-\frac{1}{(\mu_1+1)(\mu_2+1)}F_2\varphi_2(v_\mu)\\
&+\frac{1}{\mu_2+\mu_1+2}F_{1,2}\varphi_1(v_\mu) +
\frac{1}{(\mu_2+\mu_1+2)(\mu_1+1)}F_2F_1\varphi_1(v_\mu)\\
P_4(v_\mu)&=-\frac{1}{(\mu_3+1)(\mu_3+\mu_2+2)(\mu_3+\mu_2+\mu_1+3)}
\varphi_4(v_\mu)
+\frac{1}{(\mu_3+1)(\mu_2+1)(\mu_2+\mu_1+2)}F_3\varphi_3(v_\mu)\\
&-\frac{1}{(\mu_3+\mu_2+2)(\mu_1+1)}F_{2,3}\varphi_2(v_\mu)
-\frac{1}{(\mu_3+\mu_2+2)(\mu_2+1)(\mu_1+1)}F_3F_{2}
\varphi_2(v_\mu) \\
&+\frac{1}{\mu_3+\mu_2+\mu_1+3}F_{1,3}\varphi_1(v_\mu)
+\frac{1}{(\mu_3+\mu_2+\mu_1+3)(\mu_2+\mu_1+2)}F_{3}F_{1,2}\varphi_1(v_\mu)\\
&+\frac{1}{(\mu_3+\mu_2+\mu_1+3)(\mu_1+1)}F_{2,3}F_{1}\varphi_1(v_\mu)+
\frac{1}{(\mu_3+\mu_2+\mu_1+3)(\mu_2+\mu_1+2)(\mu_1+1)}F_3F_2F_1P_1v_\mu. 
\end{array}
$$
\end{remark}
\begin{theorem}\label{phicom} Let $V$ be a $\mathfrak{sl}(n+1)\ltimes 
\mathbbm{C}^{n+1}$--module,  $v_\mu\in V$ be an highest vector  of 
$\mathfrak{sl}(n+1)$ of weight $\mu$. Then    it holds
\begin{equation}
\label{phicomeq}
\varphi_i(\varphi_j(v_\mu))=\varphi_j(\varphi_i(v_\mu))\qquad 1\leq i,j\leq n+1
\end{equation}
\end{theorem}
{\bf Proof} 
Using equations (\ref{msk})  (\ref{hwveq}) and the commutation rules (\ref{comrel}) we have 
for $1\leq i\leq n$
\begin{equation}
\label{psi1i}
\begin{array}{ll}
&\varphi_{i+1}(\varphi_i(v_\mu))- \varphi_{i}(\varphi_{i+1}(v_\mu))=
-(\mu_i+2)\prod_{k=1}^{i-1}\mu_{k,i}\prod_{h=1}^{i-1}\mu_{h,i-1}P_{i+1}P_iv_\mu\\
&+\prod_{k=1}^{i-1}\mu_{k,i}\prod_{h=1}^{i-1}\mu_{h,i-1}
F_iP^2_iv_\mu+\displaystyle\sum_{\stackrel{s\leq r\leq i+1}{r+s\leq 2i}}
S^{(i+1,i)}_{r,s}P_rP_sv_\mu\\
&+(\mu_i+1)\prod_{k=1}^{i-1}\mu_{k,i}\prod_{h=1}^{i-1}\mu_{h,i-1} 
P_{i+1}P_iv_\mu
- \prod_{k=1}^{i-1}\mu_{k,i}\prod_{h=1}^{i-1}\mu_{h,i-1}
\left[P_i,F_i\right]P_iv_\mu\\
&-\prod_{k=1}^{i-1}\mu_{k,i}\prod_{h=1}^{i-1}\mu_{h,i-1}
F_iP^2_iv_\mu
+\displaystyle\sum_{\stackrel{s\leq r\leq i+1}{r+s\leq 2i}} T^{(i+1,i)}_{r,s}P_rP_sv_\mu\\
&=
\displaystyle\sum_{\stackrel{s\leq r\leq i+1}{r+s\leq 2i}}
W^{(i+1,i)}_{r,s}P_rP_sv_\mu \qquad  \mbox{ with}\ \ \ S^{(i+1,i)}_{r,s}, T^{(i+1,i)}_{r,s}, W^{(i+1,i)}_{r,s}
\in \mathcal{U}(\mathfrak{n}^-).
\end{array}
\end{equation}
 While a similar but easier computation shows for $2\leq i\leq n$ and 
$1\leq j<i-1$  that 
 \begin{equation}
\label{psi1j}
\varphi_{i+1}(\varphi_j(v_\mu))- \varphi_{j}(\varphi_{i+1}(v_\mu))=
\displaystyle\sum_{\stackrel{s\leq r\leq i+1}{r+s\leq i+j}}
W^{(i+1,j)}_{r,s}P_rP_sv_\mu \qquad   W^{(i+1,i)}_{r,s}
\in \mathcal{U}(\mathfrak{n}^-).
\end{equation} 
For $i=1$ it follows from equation (\ref{psi1i}) that
$\phi_2(\phi_1(v_\mu))= \phi_1(\phi_2(v_\mu))$ (see also Piard  \cite{Pi}).
While for $1<i\leq n$, $1\leq j< i$, using equations (\ref{invvarphi}) in the  equations (\ref{psi1i}) and (\ref{psi1j}) we obtain 
\begin{equation}
\label{Uphi}
\varphi_{i+1}(\varphi_j(v_\mu))- \varphi_{j}(\varphi_{i+1}(v_\mu))=
\displaystyle\sum_{\stackrel{s\leq r\leq i+1}{r+s\leq i+j}}
Z^{(i+1,j)}_{r,s}\varphi_r(\varphi_s(v_\mu)) \qquad   Z^{(i+1,i)}_{r,s}
\in \mathcal{U}(\mathfrak{n}^-) .
\end{equation}
Now from Theorem \ref{hwv} it follows that the vector  
$\varphi_{i+1}(\varphi_j(v_\mu))- \varphi_{j}(\varphi_{i+1}(v_\mu))$
if different from zero is a highest weight vector of $\mathfrak{sl}(n+1)$ 
with weight $\mu_{i+1,j}=\mu+\omega_{i+1}-\omega_i+\omega_j-\omega_{j-1}$,  while 
from equation  (\ref{Uphi}) follows that 
the same vector belongs to the 
$\mathfrak{sl}(n+1)$--submodule 
$$
\varphi_{i+1}(\varphi_j(v_\mu))- \varphi_{j}(\varphi_{i+1}(v_\mu))
\in \bigoplus_{\stackrel{s\leq r\leq i+1}{r+s\leq i+j}} V(\mu_{r,s})
\qquad \mu_{r,s}=\mu+\omega_{r}-\omega_{r-1}+\omega_s-\omega_{s-1}.
$$
But this is impossible because $i>1$, $s\leq r\leq i+1$, and $r+s\leq i+j$, 
imply  $\mu_{i+1,j}\neq \mu_{r,s}$.
Therefore we must have 
$$
\varphi_i(\varphi_j(v_\mu))=\varphi_j(\varphi_i(v_\mu))\qquad 1\leq i,j\leq n+1.
$$
\endpf 
 \begin{prop}\label{dirsum} Let $V$ be a 
$\mathfrak{sl}(n+1)\ltimes \mathbbm{C}^{n+1}$--module and $V(\mu)$ a 
irreducible finite $\mathfrak{sl}(n+1)$--module contained in $V$ of highest 
weight vector $v_\mu$ and highest weight $\mu$. Then for any 
there exists a $k$ such
that $\left(\varphi_{n+1}\circ \varphi_{n}\cdots \circ\varphi_1\right)^{k}(v_\mu)=0$
and the  $\mathfrak{sl}(n+1)$--modules generates by the highest weight vectors 
$\left(\varphi_{n+1}\circ \varphi_{n}\cdots \circ\varphi_1\right)^{j}(v_\mu)$,
$1\leq j< k$ are irreducible finite modules of highest weight $\mu$ such 
that for $1\leq l< j< k$ it holds 
$$
\mathcal{U}(\mathfrak{sl}(n+1))\left(\varphi_{n+1}\circ \varphi_{n}\cdots \circ\varphi_1\right)^{l}(v_\mu)\cap
\mathcal{U}(\mathfrak{sl}(n+1))\left(\varphi_{n+1}\circ \varphi_{n}\cdots \circ\varphi_1\right)^{j}(v_\mu)=
\{0\}.
$$
\end{prop}
{\bf Proof} The existence of a $k$ such that
$\left(\varphi_{n+1}\circ \varphi_{n}\cdots \circ\varphi_1\right)^{k}(v_\mu)=\{ 0\}$ and \\
$\left(\varphi_{n+1}\circ \varphi_{n}\cdots \circ\varphi_1\right)^{j}(v_\mu)
\neq \{ 0\}$ for $ 1\leq j <k$ follows from Theorem \ref{phicom} and 
the nilpotency of the operators involved. That $\left(\varphi_{n+1}\circ \varphi_{n}\cdots \circ\varphi_1\right)^{j}(v_\mu)$ is, if different from zero, 
 a highest weight vector of weight $\mu$ follows from Theorem \ref{hwv}.
 While the last statement of the Proposition holds for the same argument used in the proof of Lemma 7 of \cite{Pi}.
\endpf
From Theorem \ref{phicom}, the form of the perfect Lie algebra 
$\mathfrak{sl}(n+1)\ltimes \mathbbm{C}^{n+1}$  and the fact that from Proposition \ref{invphiP}  the action of the radical $\mathfrak{p}$ on any  highest weight vector $v_\mu$ is determined by 
the element $\varphi_i(v_\mu)$, $1\leq i\leq n+1$ we have
\begin{prop}\label{directform} Let $V$ a cyclic 
$\mathfrak{sl}(n+1)\ltimes \mathbbm{C}^{n+1}$--module with generator $v_{\mu_0}$
 then $V$ has the form
$$
V=\bigoplus_{(k_1,\dots, k_{n+1})\in \mathbbm{N}^{n+1}} 
\mathcal{U}(\mathfrak{sl}(n+1))
\varphi_{n+1}^{k_{n+1}}(\varphi_{n}^{k_{n}}\cdots(
\varphi_{1}^{k_{1}}(v_{\mu_0})\cdots ))).
$$
\end{prop} 
Since by the Weyl Theorem any cyclic 
$\mathfrak{sl}(n+1)\ltimes\mathbbm{C}^{n+1}$--module is the direct sum of finite 
dimensional irreducible modules of  $\mathfrak{sl}(n+1)$--modules we can always 
suppose that the generator $v$ belongs to an irreducible 
$\mathfrak{sl}(n+1)$--submodules. Furthermore being such module an irreducible    
$\mathfrak{sl}(n+1)$--module,  we can  suppose that $v$ is an highest 
weight vector of weight say $\mu_0$: $v=v_{\mu_0}$ as well. Such  highest 
weight vector is up multiplicative factor unique.
\begin{prop}\label{cyclicelement}
Let  $V$ be a cyclic 
$\mathfrak{sl}(n+1)\ltimes \mathbbm{C}^{n+1}$--module with generator the 
$\mathfrak{sl}(n+1)$--highest weight vector $v_{\mu_0}$ of weight $\mu_0$ then any other generators of $V$ belongs to the irreducible highest weight 
$\mathfrak{sl}(n+1)$--module 
$V(\mu_0)=\mathcal{U}(\mathfrak{sl}(n+1))v_{\mu_0}$.
\end{prop}
{\bf Proof} It is a direct consequence of Proposition  \ref{directform}
\endpf
On behalf of Corollary \ref{corrphii} and Theorem   \ref{phicom} we can give  
\begin{defi}\label{politop} Let $V$ be a cyclic  $\mathfrak{sl}(n+1)\ltimes 
\mathbbm{C}^{n+1}$--module, generated by the highest $\mathfrak{sl}(n+1)$ 
weigth vector $v_\mu\in  V$ of  weight $\mu$.\par \noindent  
Let $J^\mu_1$ be the littlest positive integer number such that
$$
\varphi^{J^\mu_1}_{1}(v_{\mu})=0.
$$
Further,  for any $0\leq i_1< J^\mu_{1}$, let   $J^\mu_{2}(i_{1})$ be  the littlest positive integer  such that
$$
 \varphi^{J_2^\mu(i_{1})}_{2}(\varphi^{i_1}_{1}
(v_\mu))=0, \quad 0\leq i_1< J^\mu_1.
$$
Recursively suppose that  we have already defined 
$J_h^\mu(i_{1},\dots i_{h-1})$ for 
$1\leq h\leq l$ and $0\leq i_r< J_{r}^\mu(i_1,\dots i_{r-1})$, $1\leq r\leq h$ 
then $J^\mu_{l+1}(i_1,i_2,\dots, i_{l})$  is the  littlest positive integer such that
$$
\varphi^{J^\mu_{l+1}(i_1,i_2,\dots, i_{l})}_{l+1}
(\varphi^{i_{l}}_{l}(\cdots (\varphi^{i_1}_1(v_\mu))))=0,
\quad 0\leq  i_1< J^\mu_1,\ \ \ 0\leq  i_{h}< J^\mu_{h}(i_{1},\dots, i_{h-1}) \quad 2\leq 
h\leq l.
$$
Hence  $J^\mu_{n+1}(i_1,i_2,\dots, i_n)$ is the littlest
positive integer such that 
$$
\varphi^{J^\mu_{n+1}(i_1,i_2,\dots, i_{n})}_{n+1}
(\varphi^{i_n}_{n}(\cdots (\varphi^{i_1}_1(v_\mu))))=0,
\quad 0\leq   i_1< J^\mu_1,\ \ \  0\leq  i_{k}< J^\mu_{k}(i_1,\dots i_{k-1}) \quad 2\leq 
k\leq n.
$$
Let us finally denote by $\mathcal{J}^\mu$ the set given by all the numbers
$J^\mu_{1}$, $J^\mu_{k}(i_1,\dots, i_{k-1})$ i.e.,:
$$
\mathcal{J}^\mu=\left\{ J^\mu_{1}, J^\mu_{k}(i_1,\dots, i_{k-1})\bigm| 
\ 0\leq  i_{1}< J^\mu_{1}, \   0\leq i_{k}< J^\mu_{k}(i_{1},\dots i_{k-1})\quad 
2\leq k\leq n+1 \right\}.
$$
\end{defi}
\begin{theorem}\label{m1mn} Let $V$ be a cyclic  $\mathfrak{sl}(n+1)\ltimes 
\mathbbm{C}^{n+1}$--module, generated by the highest $\mathfrak{sl}(n+1)$ 
weigth vector $v_\mu\in  V$ of  weight $\mu$.  Then for the elements of the set $\mathcal{J}^\mu$ of definition  \ref{politop} it holds
\begin{enumerate}
\item $J^\mu_{1}\geq 1$;
\item the set  $\mathcal{J}^\mu$,  is  bounded;
\item for $ 1\leq k\leq n$ it holds 
\begin{equation}
\label{decrescent}
J^\mu_{k}(i_{1},i_2,\dots,i_{k-1}) \leq 
J^\mu_{k}(j_1,j_2,\dots,j_{k-1})
\ \  \mbox{if }\quad   j_{l}\leq i_{l}\ \ \ 1 \leq l \leq k-1 \quad 2< k\leq n+1
\end{equation}
\item  $1\leq J^{\mu}_{k}(i_1,\dots i_{k-1})\leq  \mu_{k-1} +1\quad \forall 
i_h< J_h(i_1,\dots i_{h-1})\ \ \  1\leq h\leq k  \quad 2< k\leq n+1.$
\end{enumerate}
\end{theorem}
\hfill\eject\noindent
{\bf Proof} 
\begin{enumerate}
\item  It is obvious form the very definition of $J^{\mu}_{1}$
that it must be equal or bigger then $1$,
\item It is a consequence of  Corollary  \ref{corrphii}. 
\item   By definition 
$$
\varphi^{J^\mu_{l}(i_{1},\dots i_{l-1})-1}_{l}(\varphi^{i_{l-1}}_{l-1}(\cdots 
(\varphi^{i_k}_k(\cdots (\varphi^{i_{1}}_{1}(v_\mu))))\neq 0 \quad i_k<
J^\mu_k(i_1,\dots,i_{k-1})\ \ \ 1\leq k\leq l-1.
$$
Now from  Theorem \ref{phicom} it  follows for any $2\leq l\leq n$ and any $1\leq  k
\leq l-1$ with $i_k>0$ that 
$$
\varphi^{J^\mu_{l}(i_{1},\dots i_{l-1})-1}_{l}(\varphi^{i_{l-1}}_{l-1}(\cdots 
(\varphi^{i_k}_k(\cdot \cdot (\varphi^{i_{1}}_{1}(v_\mu))))
=\varphi_k(\varphi^{J^\mu_{l}(i_{1},\dots 
i_{l-1})-1}_{l}(\varphi^{i_{l-1}}_{l-1}(\cdots  
(\varphi^{i_k-1}_k(\cdot \cdot  (\varphi^{i_{1}}_{1}(v_\mu)))).
$$
This  implies 
$$
\varphi^{J^\mu_{l}(i_{1},\dots i_{l-1})-1}_{l}(\varphi^{i_{l-1}}_{l-1}(\cdots 
(\varphi^{i_k-1}_k(\cdots (\varphi^{i_{1}}_{1}(v_\mu))))\neq 0
$$
and in turn 
$$
J^\mu_{l}(i_{1},\cdots, i_k-1,\cdots, i_{l-1}) \geq 
J^\mu_{l}(i_{1},\cdots, i_k,\cdots, i_{l-1})\quad
 2\leq l\leq n\quad 1\leq k\leq l-1,\ \ i_k>0.
$$
\item We need only to prove  that $J^{\mu}_{k}(i_1,\dots, i_{k-1})\leq  
\mu_{k-1}+1$
for  $2\leq k\leq n+1$. Further equation (\ref{decrescent}) implies that it is 
enough to show that $J^{\mu}_{k}(0,\dots, 0)\leq  \mu_{k-1}+1$ i.e., that 
$$
\varphi_k^{\mu_{k-1}+1}(v_\mu)=0\quad 2\leq k\leq n+1.
$$
But from Theorem \ref{hwv} we 
have that if $\varphi_k^{\mu_{k-1}+1}(v_\mu)\neq 0$ then it is a highest weight vector of weight
$$
\nu=\sum_{l=1}^{k-2}\mu_l\omega_l-\omega_{k-1}
+(2\mu_{k}+1)\omega_{k}+\sum_{l=k+1}^{n+1}
\mu_l\omega_l
$$
which is impossible since $\nu$ is not a dominant integral weight of 
$\mathfrak{sl}(n+1)$.
\end{enumerate}
\endpf
\begin{defi}\label{J(V)}
According to Definition \ref{politop} and using Theorem  \ref{m1mn} and 
Proposition \ref{cyclicelement} we can associate to any cyclic $\mathfrak{sl}(n+1)\ltimes \mathbbm{C}^{n+1}$--module $V$ the   set 
\begin{equation}
\label{classet1}
\mathcal{J}(V)=\left\{\mu_0,J^{\mu_0}_{1}, J^{\mu_0}_{k}(i_{1},\dots, i_{k-1})\bigm| 
\ i_{1}< J^{\mu_0}_{1}, \   i_{k}< J^{\mu_0}_k(i_{1},\dots i_{k-1})\quad 
2\leq k\leq n+1 \right\}
\end{equation}
where $\mu_0$ is the  $\mathfrak{sl}(n+1)$ weight of the highest weight vector 
$v_{\mu_0}$ which generates $V$, and the numbers 
$J^{\mu_0}_{1}, J^{\mu_0}_{k}(i_{1},\dots, i_{k-1})$, $2\leq k\leq n+1$ 
are those defined in Definition \ref{politop}. 
\end{defi}
\begin{theorem}\label{classificationp1}
If two  cyclic finite dimensional modules $V$ and $W$ of 
$\mathfrak{sl}(n+1)\ltimes \mathbbm{C}^{n+1}$ generated respectively by the 
highest weight vector $v_{\mu_0}$ and  $w_{\nu_0}$ of weight $\mu_0$ and 
$\nu_0$ are equivalent then 
the two set $\mathcal{J}(V)$ and  $\mathcal{J}(W)$  
given by equation (\ref{classet1}) coincide. 
\end{theorem}
{\bf Proof}  
From the Propositions  \ref{directform}, and \ref{cyclicelement}  follows that the modules $V$  and $W$ decomposes as 
$$
\begin{array}{ll}
V=&V(\mu_0)\displaystyle\bigoplus_{k_l\in \mathbbm{N},\sum_{l=1}^{n+1}k_l\geq 1} \mathcal{U}(\mathfrak{sl}(n+1))
\varphi^{k_{n+1}}_{n+1}(\cdots \varphi^{k_l}_l\cdots 
(\varphi_{1}^{k_{1}}(v_{\mu_0}))))\\
W=&W(\nu_0)\displaystyle\bigoplus_{h_j\in \mathbbm{N},\sum_{j=1}^{n+1}h_j\geq 1} \mathcal{U}(\mathfrak{sl}(n+1))
\varphi^{h_{n+1}}_{n+1}(\cdots \varphi^{h_l}_l\cdots 
(\varphi_{1}^{h_{1}}(w_{\nu_0}))))
\end{array}
$$
here  $V(\mu_0)$ (resp. $W(\nu_0)$) is the irreducible $\mathfrak{sl}(n+1)$--module 
with highest weight $\mu_0$ (res. $\nu_0$) and  highest weight vector $v_{\mu_0}$ (resp. $w_{\nu_0}$).\par 
If the module $V$ and $W$ are equivalent then there exists an invertible linear
operator $T\in \mbox{Hom}_{\mathbbm{C}}(V,W)$ such that $XT=TX$ for any element $X$ in $\mathcal{U}(\mathfrak{sl}(n+1)\ltimes \mathbbm{C}^{n+1})$. Hence, 
since  $V(\mu_0)=\mathcal{U}(\mathfrak{sl}(n+1))v_{\mu_0}$ and 
$W(\nu_0)=\mathcal{U}(\mathfrak{sl}(n+1))w_{\nu_0}$ the Schur Lemma \cite{H} 
\cite{Kn}  
implies
 that $\mu_0=\nu_0$ and therefore that $W(\nu_0)\simeq V(\mu_0)$. 
For any highest weight vector $v_\mu$ in $V$,  $Tv_\mu$ is an highest weight 
vector in $W$ with the same weight $\mu$. Hence, since T commutes with the action of $\mathfrak{sl}(n+1)\ltimes \mathbbm{C}^{n+1}$ we have 
$T\varphi_i(u_\mu)=\varphi_i(Tu_\mu)$ for every highest weight vector 
$u_\mu$ and any $i$, $1\leq i\leq n+1$. This implies 
$$
\varphi^{k_{n+1}}_{n+1}(\cdots \varphi^{k_l}_l\cdots 
(\varphi_{1}^{k_{1}}(v_{\mu_0}))))\neq 0 \Longleftrightarrow 
 \varphi^{k_{n+1}}_{n+1}(\cdots \varphi^{k_l}_l\cdots 
(\varphi_{1}^{k_{1}}(w_{\mu_0})))\neq 0.
$$
This latter equation obviously shows that  set $\mathcal{J}(V)$ and 
$\mathcal{J}(W)$ coincide.\par 
\endpf
In order to complete the classification of all cyclic 
$\mathfrak{sl}(n+1)\ltimes \mathbbm{C}^{n+1}$--modules,
it remains to show that for any set $\mathcal{M}$  which satisfies the 
requirements of Definition \ref{J(V)} together with those of Theorem \ref{m1mn}
there exists a cyclic $\mathfrak{sl}(n+1)\ltimes \mathbbm{C}^{n+1}$--module
$V$ such that $\mathcal{J}(V)=\mathcal{M}$.\par 
This will be done in the next section where we shall show how 
such modules are quotient modules of the restriction to 
$\mathfrak{sl}(n+1)\ltimes \mathbbm{C}^{n+1}$ of finite dimensional irreducible
$\mathfrak{sl}(n+2)$--modules.
\section{The $\mathfrak{sl}(n+2)$--modules as  
$\mathfrak{sl}(n+1)\ltimes \mathbbm{C}^{n+1}$--modules} 
Viewed in the light of the previous section the embeddings  $\Phi$ and 
$\Theta$ of
$\mathfrak{sl}(n+1)\ltimes \mathbbm{C}^{n+1}$ in 
$\mathfrak{sl}(n+2)$ given in Theorem \ref{authomorphism} by the formulas 
(\ref{emmb1}) and (\ref{emmb2}) are 
$\mathfrak{sl}(n+1)\ltimes \mathbbm{C}^{n+1}$--modules $\mathbbm{C}^{n+2}_\Phi$ and
$\mathbbm{C}^{n+2}_\Theta$ 
whose set $\mathcal{J}(\mathbbm{C}^{n+2}_\Phi)$ and 
$\mathcal{J}(\mathbbm{C}^{n+2}_\Theta)$ are respectively
$$
\begin{array}{ll}
&\mathcal{J}(\mathbbm{C}^{n+2}_\Phi)=\left\{\mu_0=0,\ \ J_{1}=2 \ \ \ 
J^{\mu_0}_{k+1}(\underbrace{0,\dots,0}_{k-times})= J^{\mu_0}_{k+1}(\underbrace{1,\dots,0}_{k-times})=1, 
\ \ 1\leq k\leq n\right\}\\
&\\
&\mathcal{J}(\mathbbm{C}^{n+2}_\Theta)
=\left\{\mu_0=\omega_{n},\ \ J^{\mu_0}_{1}=1, \ \  
J^{\mu_0}_{k+1}(\underbrace{0,\dots,0}_{k-times})=1,\ \ 
1\leq k< n,\  J^{\mu_0}_{n+1}(\underbrace{0,\dots,0}_{n-times})=2\right\}.
\end{array}
$$
The aim of this section is to investigate the restrictions through such 
embeddings of the finite dimensional irreducible 
$\mathfrak{sl}(n+2)$--modules to 
$\mathfrak{sl}(n+1)\ltimes \mathbb{C}^{n+1}$.\par 
Observe that the  automorphism $\Xi$ defined in Proposition \ref{propemb} 
can be extended to the whole weight space $\mathfrak{h}^*$ by setting
$$
\Xi(\lambda)=\Xi(\sum_{k=1}^{n+1}\lambda_k\omega_k)= 
\sum_{k=1}^{n+1}\lambda_{n+2-k}\omega_k,\qquad \lambda \in \mathfrak{h}^*
$$
 and that therefore
it holds
\begin{prop}\label{equi}
Let $V(\lambda)$ be  a irreducible finite dimensional 
$\mathfrak{sl}(n+2)$--module, then $V(\lambda)$ viewed as  
$\mathfrak{sl}(n+1)\ltimes\mathbbm{C}^{n+1}$--module using the embedding $\Phi$ 
(\ref{emmb1}) is equivalent  to the module $V(\Xi(\lambda))$ always viewed as  
$\mathfrak{sl}(n+1)\ltimes\mathbbm{C}^{n+1}$--module using 
the embedding $\Theta$ (\ref{emmb1})
\end{prop} 
Therefore, since the automorphism $\Xi$ preserves the integral dominant weight 
of  $\mathfrak{sl}(n+2)$, in order  to obtain all the   modules of  
$\mathfrak{sl}(n+1)\ltimes \mathbbm{C}^{n+1}$, given   by the restriction to 
it of the irreducible finite dimensional $\mathfrak{sl}(n+2)$--modules it 
is enough   to consider one of the two 
embedding, for our convenience let us chose $\Phi$ (\ref{emmb1}).\par
Although this fact would allow us to simply talk about restriction of 
$\mathfrak{sl}(n+2)$--modules to  $\mathfrak{sl}(n+1)\ltimes \mathbbm{C}^{n+1}$
in order to avoid confusion we shall denote by $V(\lambda)_\Phi$ the 
restriction to $\mathfrak{sl}(n+1)\ltimes \mathbbm{C}^{n+1}$ of 
the irreducible $\mathfrak{sl}(n+2)$--module $V(\lambda)$.\par 
Observe that if $u_\mu$ is a weight vector of $\mathfrak{sl}(n+2)$ on $V(\lambda)$ then it is a weight vector  of $\mathfrak{sl}(n+1)$ in 
$V(\lambda)_\Phi$ of weight $\mu=\sum_{k=2}^{n+1}\mu_k\omega_k$. In particular if   $u_\mu$ is a highest weight vector of $\mathfrak{sl}(n+1)$ in 
$V(\lambda)_\Phi$  
then using the embedding (\ref{emmb1}) the element $\varphi_i(u_\mu)$, $1\leq i\leq n+1$ (\ref{hwveq}) are 
\begin{equation}
\label{carphivlambda}
\varphi_{i}(u_\mu)=F_{1,i}^{((\mu(i))}(u_\mu)
+\sum_{k=1}^{i-1}Q^{(\mu(i))}_{k+1}F^{((\mu(i))}_{1,k}(u_\mu)\qquad
1\leq i\leq n+1 
\end{equation}
where the elements $Q^{(\mu(i))}_{k+1}$, $1\leq k\leq n=1$  are those given by  Definition \ref{Qp} applied to the universal enveloping algebra $\mathcal{U}(\mathfrak{sl}(n+2))$. For instance the first elements $\varphi_i(u_\mu)$ become 
$$
\begin{array}{ll}
\varphi_1(u_\mu)&=F_{1} u_\mu\\
\varphi_2(u_\mu)&=-(\mu_2+1)F_{1,2}u_\mu+F_2F_{1}u_\mu\\
\varphi_3(u_\mu)&=(\mu_3+1)(\mu_3+\mu_2+2)F_{1,3}
 u_\mu-(\mu_3+\mu_2+2)F_3F_{1,2}u_\mu\\
&-(\mu_3+1)F_{2,3}F_1u_\mu +F_3F_2F_{1}u_\mu\\
\varphi_4(u_\mu)&=-(\mu_4+1)(\mu_4+\mu_3+2)(\mu_4+\mu_3+\mu_2+3)F_{1,4}
u_\mu\\
&+(\mu_4+\mu_3+2)(\mu_4+\mu_3+\mu_2+3)F_4F_{1,3}u_\mu\\
&+(\mu_4+1)(\mu_4+\mu_3+\mu_2+3)F_{3,4}F_{1,2}u_\mu 
-(\mu_4+\mu_3+\mu_2+3)F_{4}F_{3}F_{1,2}u_\mu\\
&+(\mu_4+1)(\mu_4+\mu_3+2)F_{2,4}F_{1}u_\mu
-(\mu_4+\mu_3+2)F_4F_{2,3}F_1u_\mu
\\
& 
-(\mu_4+1)F_{3,4}F_2F_1u_\mu+ 
F_4F_3F_2F_1u_\mu. 
\end{array}
$$
\begin{theorem}\label{cyclicimb} For any integral dominant weight 
$\lambda=\sum_{i=1}^{n+1}\lambda_i\omega_i$ 
of  $\mathfrak{sl}(n+2)$ the restriction $V(\lambda)_\Phi$ of  the irreducible 
$\mathfrak{sl}(n+2)$--module 
$V(\lambda)$ to $\mathfrak{sl}(n+1)\ltimes \mathbb{C}^{n+1}$ is a cyclic module, 
with generator $v_\lambda$ the $\mathfrak{sl}(n+2)$--highest weight vector in 
$V(\lambda)$,  whose   $\mathfrak{sl}(n+1)$--weight is \\
$\lambda_0=\sum_{k=2}^{n+1}\lambda_k\omega_k$. \par 
\end{theorem}
{\bf Proof} 
Define on the set $\Delta^+$  of positive roots of $\mathfrak{sl}(n+2)$  a total order $\succ_n$:
\beq
\label{order}
\alpha_{p,q}\succ_n \alpha_{r,s} \qquad \Longleftrightarrow \quad \mbox{if}\ \ \ \ p< r \quad \mbox{or} \ \ \ \ p=r \ \ \mbox{and} \ \  q< s
 \eeq
Correspondingly, we say that $F_{p,q}\succ_n F_{r,s}$ if $\alpha_{p,q}\succ_n\alpha_{r,s}$ so
$$
F_{1,1}\succ_n \cdots F_{1,n+1}\succ_n F_{2,2}\succ_n F_{2,3}\succ_n\cdots \succ F_{n-1,n+1}\succ_n F_{n,n}\succ_n  F_{n,n+1}\succ_n F_{n+1,n+1}.
$$
Then using this ordering the  $\mathfrak{sl}(n+2)$--module $V(\lambda)$ 
is the span  of vectors of the type \cite{H} \cite{FFL}   
$$
F^{a_{n+1}}_{n+1}F^{a_{n,n+1}}_{n,n+1}\cdots F^{a_{i,j}}_{i,j}\cdots  
F^{a_{1,n}}_{1,n}F^{a_{1,n-1}}_{1,n-1}\cdots F^{a_1}_1v_\lambda 
\quad a_{i,j}\geq 0.
$$
Applying the embedding $\Phi$ (\ref{emmb1}) this means that 
$$
V(\lambda)=\mathcal{U}(\mathfrak{sl}(n+1))\mathcal{U}(\mathfrak{p}) v_\lambda
$$
i.e., that $V(\lambda)$ is a  cyclic 
$\mathfrak{sl}(n+1)\ltimes \mathbbm{C}^{n+1}$--module with generators $v_\lambda$.\par \noindent 
\endpf
\begin{theorem}\label{Jlambda} Let $V(\lambda)$ be an  irreducible finite dimensional $\mathfrak{sl}(n+2)$--module of highest weight 
$\lambda=\sum_{k=1}^{n+1} \lambda_k\omega_k$. Then for the corresponding 
cyclic module $V(\lambda))_\Phi$ of  
$\mathfrak{sl}(n+1)\ltimes \mathbbm{C}^{n+1}$ we have 
\begin{equation}
\label{jvlambda}
\mathcal{J}(V(\lambda)_\Phi))=\left\{
\begin{array}{ll}
&\lambda_0=\sum_{i=2}^{n+1} \lambda_i\omega_i,\ \ 
J^{\lambda_0}_{1}=\lambda_{1}+1,\\
&J^{\lambda_0}_k(i_{1},\cdots i_{k-1})=\lambda_{k}+1
\end{array}
\vline 
\begin{array}{ll} 
&0\leq i_h\leq \lambda_{h}, \    1\leq h\leq n+1\\
& 2\leq  k\leq n+1
\end{array}
\right\}. 
\end{equation}
\end{theorem}
{\bf Proof}
Let $v_{\lambda_0}$ be the highest weight vector of $\mathfrak{sl}(n+2)$--module $V(\lambda)$ seen 
as  highest weight vector of $\mathfrak{sl}(n+1)$ with weight 
$\lambda_0$.\par 
 Since from equations (\ref{emmb1}) and (\ref{hwveq}) we have
$$
\varphi^k_1(v_{\lambda_0})=F^k_1(v_{\lambda_0}),
$$
formula $J^{\lambda_0}_1=\lambda_1+1$ follows from the theory of finite dimensional   
$\mathfrak{sl}(2)$--modules.\par 
We can now proceed by induction suppose that we have already shown \\
$J^{\lambda_0}_h(i_{1},\cdots i_{h-1})= \lambda_h+1$ for $h\leq k$, we want show that 
 $J^{\lambda_0}_{k+1}(i_{1},\cdots, i_{k})= \lambda_{k+1}+1$. 
First we shall show that 
$$
\varphi^{i_{k+1}}_{k+1}(\varphi^{i_{k}}_{k}(\cdots 
(\varphi^{i_1}_1(v_{\lambda_0})))\neq 0 \quad 
\mbox{if $i_{k+1}\leq \lambda_{k+1}$, 
$0\leq i_h\leq \lambda_{h}$, $1\leq h\leq k$.}
$$
Let $u_\eta$ be  defined by
$$
u_\eta=\varphi^{i_{k}}_{k}(\cdots 
(\varphi^{i_1}_1(v_{\lambda_0})). 
$$
We want show that if $\lambda_{k+1}\geq 1$ then 
$$
\varphi_{k+1}(u_\eta)\neq 0.
$$
From Theorem \ref{hwv} the vector $u_\eta$ is a highest weight vector of 
$\mathfrak{sl}(n+1)$ of $\mathfrak{sl}(n+2)$ weight 
\begin{equation}
\label{muhw}
\eta=\sum_{j=1}^{n+1}\eta_k\omega_k=
\sum_{j=1}^{k-1}(\lambda_j+i_{j}-i_{j+1})\omega_j+(\lambda_{k}+i_k)\omega_{k}
+\sum_{j=k+1}^{n+1}\lambda_j\omega_j,
\end{equation}
and therefore of $\mathfrak{sl}(n+1)$ weight 
$\eta^{(2)}=\sum_{j=2}^{n+1}\eta_k\omega_k$.
Hence   using 
equations (\ref{emmb1}) and (\ref{carphivlambda}) we have
$$
\varphi_{k+1}(u_\eta)=F_{1,k+1}^{(\eta(k+1))}u_\eta
+\sum_{j=1}^{k}Q^{(\eta(k+1))}_{j+1}F^{(\eta(k+1))}_{1,j}u_\eta.
$$
Since $\eta_{k+1}=\lambda_{k+1}\geq 1$, from the theory of the finite dimensional $\mathfrak{sl}(2)$--module  (and from Theorem \ref{dyckbasis} as 
well) follows that $F_{1,k+1}^{(\eta(k+1))}u_\eta\neq 0$.  Hence  if 
$Q^{(\eta(k+1))}_{j+1}F^{(\eta(k+1))}_{1,j}u_\eta=0$ for any $1\leq j\leq k$ 
we have  
$$
\varphi_{k+1}(u_\eta)=F_{1,k+1}^{(\eta(k+1))}u_\eta\neq 0.
$$
We may suppose, therefore, without loosing generality, that $Q^{(\eta(k+1))}_{j+1}F^{(\eta(k+1))}_{1,j}u_\eta\neq 0$ for $1\leq j\leq k$. 
From Proposition \ref{Qp} and formula (\ref{emmb1}) we have that  
 \begin{equation}
\label{Qj+1}
Q^{(\eta(k+1))}_{j+1}=R^{(\eta(k+1))}_{j+1,k+1-j}=
\sum_{l=0}^{k-j}F^{(\mu(k+1))}_{j+1+l,k+1}R^{(\eta(k+1))}_{j+1,l}
\end{equation}
and therefore
\begin{equation}
\label{qf1}
 Q^{(\eta(k+1))}_{j+1}F^{(\eta(k+1))}_{1,j}u_\eta =
\sum_{l=0}^{k-j}F^{(\mu(k+1))}_{j+l+1,k+1}R^{(\eta(k+1))}_{j+1,l}
F^{(\eta(k+1))}_{1,j}u_\eta.
\end{equation}
Consider now the  simple subalgebra of $\mathfrak{sl}(n+2)$ generated by the elements
$E_{p,q}$,$F_{p,q}$, with $1\leq p\leq q\leq k$, which is  isomorphic to 
$\mathfrak{sl}(k+1)$. The vector $u_\eta$ is also a highest weight vector of
such simple subalgebra with highest weight vector 
$\eta^{(k)}=\sum_{l=1}^k\eta_l\omega_l$.\par 
Now form Definition \ref{rr}  we have that the vectors 
$R^{(\eta(k+1))}_{j+1,l}F^{(\eta(k+1))}_{1,j}u_\eta$  belong to the 
the $\mathfrak{sl}(k+1)$-module generated by $u_\eta$. Hence we can write them  as linear combination of the element of the (FFL) basis \ref{dyckbasis}  of 
the irreducible finite dimensional $\mathfrak{sl}(k+1)$--module 
$V(\eta_{k})$  of highest weight $\eta_k$ and highest weight vector $u_\eta$:
$$
R^{(\eta(k))}_{j+1,l}F^{(\eta(k+1))}_{1,j}u_\eta=
\sum_{\mathbf{s}_{k}\in S(\eta_k)} c_{\mathbf{s}_{k}}
F^{\mathbf{s}_k} u_\eta
$$
where  $S(\eta_k)$ is the set defined in Theorem  \ref{dyckbasis}. 
Substituting this last equation in (\ref{qf1}) yields 
$$
 Q^{(\eta(k+1))}_{j+1}F^{(\eta(k+1))}_{1,j}u_\eta =
\sum_{l=0}^{k-j}\sum_{\mathbf{s}_{k}\in S(\eta_k)} F^{(\eta(k+1))}_{j+l+1,k+1} F^{\mathbf{s}_k} u_\eta= \sum_{l=j+1}^{k+1}\sum_{\mathbf{s}_{k}\in S(\eta_k)} F^{(\eta(k+1))}_{l,k+1} F^{\mathbf{s}_k} u_\eta.
$$
 We claim that if $F^{\mathbf{s}_k}u_\eta$
is an element of the Feigin Fourier Littelmann basis of the irreducible finite 
dimensional $\mathfrak{sl}(k+1)$--module $V(\eta_{k})$  then
 $F^{\eta(k+1)}_{l,k+1} F^{\mathbf{s}_k} u_\eta$ 
 is element of the 
Feigin Fourier Littelmann basis of the $\mathfrak{sl}(n+2)$--module 
$V(\eta)$  of highest weight $\eta$.
Let $\mathbf{s}=(\mathbf{s}_{k},\mathbf{s}_{l,k+1})$  the multi--exponent
of $\mathfrak{sl}(n+2)$ such that 
$F^{\mathbf{s}}u_\eta=F^{(\eta(k+1))}_{l,k+1} F^{\mathbf{s}_k} u_\eta$
for Theorem \ref{dyckbasis} we have to show that for any 
Dyck path $\mathbf{p}=(\beta(0),\dots, \beta(h))$ with say  $\beta(0)=\alpha_r$  and $\beta(h)=\alpha_s$ it holds 
\begin{equation}
\label{dyckformula1}
s_{\beta(0)}+ s_{\beta(1)}+\cdots 
+s_{\beta(h)}\leq \sum_{i=r}^{s}\eta_i.
\end{equation}
Let $\Delta^+_{i}$, $1\leq i\leq n+1$  be  the subset of the set of the positive 
roots 
$\Delta^+$ of $\mathfrak{sl}(n+2)$ which can be written as linear combination of 
the simple roots $\alpha_l$ with $1\leq l\leq i$ (i,e,  the set of the 
positive roots of the subalgebra $\mathfrak{sl}(i+1)$).   
Then for any Dyck path $\mathbf{p}=(\beta(0),\dots, \beta(h))$ of
$\mathfrak{sl}(n+2)$ we have 
\begin{equation}
\label{sprop}
(s_{k})_{\beta(i)}=0 \mbox{ if }\  \beta(i) \notin \Delta^+_{k}
\quad \mbox{and}  \quad 
 (s_{l,k+1})_{\beta(i)}=0 \mbox{ if }\ \beta(i)\notin \Delta^+_{k+1}. 
\end{equation}
Now  if the elements of the Dyck path $\mathbf{p}$ belong to $\Delta^+_k$ 
i.e, if $\beta(h)=\alpha_s$ with $1\leq s\leq k$ equation (\ref{dyckformula}) 
is verified because $\mathbf{s}_k\in S(\eta_k)$. While, if this is not the case, 
using equation (\ref{sprop}),  we may assume that $\beta(h)=\alpha_{k+1}$, then  equation (\ref{dyckformula1}) becomes $\sum_{i=r}^{k+1} 
s_{\beta(i)}\leq \sum_{i=r}^{k+1}\eta_i$.
Since $\eta_{k+1}\geq 1$ if $\beta(0)=\beta(h)=\alpha_{k+1}$ 
equation (\ref{dyckformula1}) is obviously verified since $(s_{l,k+1})_{\alpha_{k+1}}\leq 1$
otherwise let be $(\beta(0),\dots,\beta(l))$ the subset given by the root of 
$\mathbf{p}$ which belong to $\Delta^+_k$ 
 then from the very definition \ref{Dyck} of Dyck path it follows that there exist a $p$ with $1\leq p\leq k$ such that 
$\beta(l)=\alpha_{p,k}$.  Then  adding to the end of the sequence
 $(\beta(0),\dots,\beta(l))$ 
the elements $(\alpha_{p+1,k},\alpha_{p+2,k}, \dots ,\alpha_{k,k})$,
we make it a Dyck path of $\mathfrak{sl}(k+1)$.
Hence we have that equation (\ref{dyckformula1}) is satisfied because 
$$
\sum_{j=0}^ls_{\beta(j)}
+\sum_{j=l+1}^{h}s_{\beta(j)}\leq 
\sum_{j=0}^l(s_k)_{\beta(j)}
+(s_{1,k+1})_{\alpha_{l,k+1}}
\leq 
\sum_{j=0}^l(s_k)_{\beta(j)}+ 
\sum_{t=1}^{k-p}(s_k)_{\alpha_{p+t,k}}+1\leq 
\sum_{i=r}^k\eta_k +\eta_{k+1}
$$
 being  $\mathbf{s}_k\in S(\eta_k)$ and $\eta_{k+1}\geq 1$.\par  
Now  $\eta_{k+1}>1$ and Theorem \ref{dyckbasis} imply also that 
$F_{1,k+1}^{(\eta(k+1))}u_\eta$ is  an element of the (FFL) 
basis of the $\mathfrak{sl}(n+2)$--module   $V(\eta)$.\par 
Therefore $\varphi_{k+1}(u_\eta)$ can be write  
as linear combination  of elements such basis as 
\begin{equation}
\label{phiuk}
\varphi_{k+1}(u_\eta)=F_{1,k+1}^{(\eta(k+1))}u_\eta+
\sum_{j=1}^{k}\sum_{l=j+1}^{k+1}\sum_{\mathbf{s}_{k}\in S(\eta_k)} F^{\eta(k+1)}_{l,k+1} F^{\mathbf{s}_k} u_\eta.
\end{equation}
Finally, since the element $F_{1,k+1}^{(\mu(k+1))}u_\eta$  of the (FFL) basis
 appears only once in (\ref{phiuk}) and it is different from zero we have 
$\varphi_{k+1}(u_\eta)\neq 0$.\par 
Summarizing we have shown that 
$$
\lambda_{k+1}\geq 1\Longrightarrow F_{1,k+1}^{((\eta(k+1))}(u_\eta)\neq 0\Longrightarrow 
\varphi_{k+1}(u_\eta)\neq 0.
$$
Set now $u^l_\eta=\varphi^l_{k+1}(u_\eta)$, $0\leq l\leq \lambda_{k+1}$ 
($u^0_\eta=u_\eta$) the same 
argument shows that 
$$
\lambda_{k+1}\geq l\Longrightarrow   F_{1,k+1}^{((\eta(k+1))}(u^{l-1}_\eta)\neq 0\Longrightarrow 
\varphi_{k+1}(u^{l-1}_\eta)\neq 0.
$$
Therefore for any  $i_h\leq \lambda_h$,  $1\leq h\leq k$ we have
$J^{\lambda_0}_{k+1}(i_1,\dots, i_{k})\geq \lambda_{k+1}+1$.\par 
For the Theorem \ref{m1mn} (3), to prove that $J^{\lambda_0}_{k+1}(i_1,\dots, i_{k})\leq  
\lambda_{k+1}+1$ we need only to show that
$$
\varphi^{\lambda_{k+1}+1}(v_{\lambda_0})=0.
$$
But arguing like  in the proof of Theorem \ref{m1mn} we see that if 
$\varphi^{\lambda_k+1}(v_{\lambda_0})$ is different from zero then there  
would be 
a highest weight vector of $\mathfrak{sl}(n+1)$ with a non dominant weight, 
which is impossible.
\endpf
\begin{cor}\label{dec}
 The  cyclic  $\mathfrak{sl}(n+1)\ltimes\mathbbm{C}^{n+1}$--module $V(\lambda)_\Phi$   decomposes as  $\mathfrak{sl}(n+1)$--module as 
 $$
V(\lambda)_\Phi=\bigoplus_{l=1}^{n+1}\bigoplus_{k_l=1}^{\lambda_l}
V\left(\lambda_0+\sum_{l=1}^n(k_l-k_{l+1})\omega_l\right).
$$
\end{cor}
\begin{defi}\label{setM}
Let $\mathfrak{M}$ be the collection of all sets $\mathcal{M}=(\mu_0,\mathcal{M}^{\mu_0})$ such that 
 $\mu_0=\sum_{l=1}^n\mu_k\omega_k$ is a dominant weight of $\mathfrak{sl}(n+1)$ and   $\mathcal{M}^{\mu_0}$ is a set of positive integers defined as 
follows:\par\noindent
\begin{equation}
\label{setMeq}
\mathcal{M}^{\mu_0}=\left\{
\begin{array}{ll}
M_1^{\mu_0}, M_k^{\mu_0}(i_1,\dots i_{k-1}),&\\
\ 2\leq k\leq n+1&
\end{array}\vline
\begin{array}{ll} 
0\leq i_k < M^{\mu_0}_k(i_1,\dots ,i_{k-1}), \ \ 1\leq k\leq n+1;&\\
M_1^{\mu_0}\geq 1;&\\
1\leq   M_k^{\mu_0}(i_1,\dots i_{k-1})\leq \mu_{k-1}+1,\ \ 2\leq k\leq n+1;
&\\
 M_k^{\mu_0}(i_1,\dots,i_{k-1})\geq M_k^{\mu_0}(j_1,\dots, j_{k-1})&\\
\mbox {if} \ \  i_l\leq j_l,\ \ 0\leq l\leq k-1 \ \  2\leq k\leq n+1.&
\end{array}\right\}
\end{equation}
\end{defi}
\begin{theorem}\label{setjlambda} For any set $\mathcal{M}\in \mathfrak{M}$  there exist a cyclic module $V$ of
$\mathfrak{sl}(n+1)\ltimes \mathbbm{C}^{n+1}$ such that
$$
\mathcal{M}=\mathcal{J}(V).
$$
Such  $\mathfrak{sl}(n+1)\ltimes \mathbbm{C}^{n+1}$--module $V$ can be constructed as 
quotient of a suitable \\ $\mathfrak{sl}(n+1)\ltimes \mathbbm{C}^{n+1}$--module
 $V(\lambda)_\Phi$ obtained as restriction to  
$\mathfrak{sl}(n+1)\ltimes \mathbbm{C}^{n+1}$ of a finite dimensional irreducible 
 $\mathfrak{sl}(n+2)$--module $V(\lambda)$.
\end{theorem}
{\bf Proof} Let $\mathcal{M}$ be a set in $\mathfrak{M}$  (see  Definition \ref{setjlambda})
and let  $\lambda_0=\sum_{k=1}^n\lambda^0_k\omega_k$ be the integral dominant weight in $\mathcal{M}$.\par 
Define the dominant integral weight  $\lambda$ of  $\mathfrak{sl}(n+2)$ as follows
$$
\lambda=\sum_{k=1}^{n+1}\lambda_k\omega_k\quad \lambda_1=M^{\lambda_0}_1-1, \ \lambda_{k+1}=\lambda^0_{k-1}\ \  2\leq k\leq n
$$
where $M^{\lambda_0}_1$ is the positive integer which appears in the definition of $\mathcal{M}$. Let $V(\lambda)$ be  the associated irreducible finite dimensional  $\mathfrak{sl}(n +2)$--module with highest weight vector $v_\lambda$,  and let $V(\lambda)_\Phi$ be its restriction 
to  $\mathfrak{sl}(n+1)\ltimes \mathbbm{C}^{n+1}$.
From equation (\ref{jvlambda}) and the Definition  \ref{setM}  follows that
$$
\begin{array}{ll}
&M^{\lambda_0}_1=\lambda_1+1\\
&M^{\lambda_0}_k(i_1,\dots,i_{k-1})\leq \lambda_k+1 \quad 0\leq i_h < M^{\mu_0}_h(i_1,\dots ,i_{h-1}),\ 1\leq h\leq k-1
\quad 2\leq k\leq n+1
\end{array}
$$
For any integer $k$ with $1\leq k\leq n$, let  $\mathcal{I}_k=\left\{i_1,\dots, i_k\right\}$ be the set of positive integer numbers such that 
$i_l<M^{\lambda_0}_l(i_1,\dots,i_{l-1})$, $1\leq l\leq k$  
and $M^{\lambda_0}_{k+1}(i_1,\dots,i_{k})$, is strictly less then  
$\lambda_{k}+1$.   Further for any $1\leq k\leq n+1$ and any $k$-uple  
$\left\{i_1,\dots, i_k\right\}\in \mathcal{I}_k$
let  $u^{(k)}(i_1,\dots. i_{k})$ be 
defined by
$$
u^{(k)}(i_1,\dots,i_{k})
=\varphi^{M^{\lambda_0}_{k+1}(i_1,\dots,i_{k-1})+1}(\varphi^{i_{k-1}}_{i_{k-1}}(\cdots 
(\varphi^{i_1}_1(v_\lambda)))).
$$
Finally  for any $1\leq k\leq n+1$ and any $k$-uple  
$\left\{i_1,\dots, i_k\right\}\in \mathcal{I}_k$  
denote by $W^{(k)}(i_1,\dots,i_{k})$ the 
$\mathfrak{sl}(n+1)\ltimes \mathbbm{C}^{n+1}$--submodule  generated by 
$u^{(k)}(i_1,\dots, i_{k})$:
$$
W^{(k)}(i_1,\dots,i_{k})=
\mathcal{U}(\mathfrak{sl}(n+1)\ltimes \mathbbm{C}^{n+1})u^{(k)}(i_1,\dots i_{k})
$$
and with $W_\mathcal{M}$ their union
$$
W_\mathcal{M}=\bigcup_{k=1}^{n}\left(\bigcup_{(i_1,\dots, i_k)\in \mathcal{I}_k}W^{(k)}(i_1,\dots,i_{k})\right).
$$
Let $V_\mathcal{M}$ be finally the quotient module
$$
V_\mathcal{M}=(V(\lambda))_\Phi/W_\mathcal{M}.
$$
Then by construction from  Theorems \ref{m1mn},   \ref{jvlambda} and the 
Definition \ref{setM} it follows 
$$
\mathcal{J}(V_\mathcal{M})=\mathcal{M}.
$$
\endpf
Putting together  Theorem \ref{classificationp1} and Theorem \ref{setjlambda} we have
\begin{theorem}\label{classificationth} 
The class of the finite dimensional 
$\mathfrak{sl}(n+1)\ltimes \mathbbm{C}^{n+1}$--modules is in one to one corresponence with the collection of sets $\mathfrak{M}$ defined in Definition \ref{setM}.
\end{theorem}

\end{document}